\documentclass[a4,10pt]{article}
\usepackage{calc}
\usepackage[all]{xy}
\usepackage[centertags]{amsmath}
\usepackage{latexsym}
\usepackage{amsfonts}
\usepackage{amssymb}
\usepackage{amsthm}
\usepackage[all]{xy}
\usepackage{fancyhdr}
\usepackage [dvips]{epsfig}

\usepackage{newlfont}
\usepackage[latin1]{inputenc}
\usepackage[french,english]{babel}
\usepackage{graphicx}
\usepackage{t1enc}
\theoremstyle{plain}
\newtheorem{lem}{Lemma}[section]
\newtheorem{prop}{Proposition}[section]
\newtheorem{theo}{Theorem}[section]
\newtheorem{definition}{Definition}[section]

\newtheorem{rem}{Remark}[section]

\theoremstyle{plain}
\begin{document}
\title{Poisson (co)homology of polynomial Poisson algebras in \\
dimension four : Sklyanin's case}
\author{Serge Rom\'eo Tagne Pelap}
\maketitle
\begin{abstract}
In this paper, we compute the Poisson (co)homology of  a polynomial Poisson structure given by two Casimir polynomial functions which define a complete intersection with an isolated singularity.\\
\end{abstract}
\nocite{*}
\section*{Introduction}
The canonical  or Poisson homology was introduced independently by Brylinsky (~\cite{bry}) (as an important tool in computations of Hochschild and cyclic homology), and by Koszul and Gelfand-Dorfman (inspired by their algebraic approach to the study of bi-hamiltonian structures).\\
This homology is defined as the homology of a differential complex degree -1  $\{\Omega^\bullet(M), \partial_{\pi}\}$ on a Poisson manifold $(M, \pi)$, where $\pi$ is a Poisson structure given either by a $2$-tensor field $\pi\in H^0(M, \wedge^2TM)$ or, by the corresponding Poisson bracket on algebra of functions (smooth, algebraic,...) of $M$ : $\{f, g\}=<df\wedge dg, \pi>.$ The $(-1)-$differential $\partial_{\pi}$ is given on the decomposable differential forms as
$$\partial_k(f_0df_1\wedge...\wedge df_k)=\displaystyle{\sum_{1\leq i\leq k}}(-1)^{i+1}\{f_0, f_i\}df_1\wedge...\wedge \widehat{df_i}\wedge ...\wedge df_k$$
$$+\displaystyle{\sum_{1\leq i<j\leq k}} (-1)^{i+j}f_0d\{f_i,f_j\}\wedge df_1\wedge...\wedge \widehat{df_i}\wedge...\widehat{df_j}\wedge...\wedge df_k$$
and acts from $\Omega^k(M)$ to $\Omega^{k-1}(M).$\\
There is a dual notion of Poisson cohomology introduced by Lichnerowicz (~\cite{lic}). But this duality is quite subtle : if the Poisson homology is non-degenerate (symplectic) then the Poisson homology and Poisson cohomology admit a Poincaré-duality ismorphism (~\cite{bry}) which is extended to a class of so-called unimodular Poisson manifolds (those whose Weinstein modular cohomology class is trivial (~\cite{xu1}) ).\\
In general, the Poisson homology has very bad functorial properties (this was observed in many papers (~\cite{gin}) ) which make its computation a challenging and difficult problem. The Poisson (co)homology depends heavily on the properties of the initial Poisson structure and, basically, on the structure of its degeneration locus. "Simple" Poisson structures (like symplectic ones or the slightly more difficult case of regular structures which have a constant rank symplectic foliation) have well-studied Poisson (co)homology. \\
An interesting and difficult problem is to compute the corresponding (co)homology groups in the case of "quadratic" Poisson structures.\\
The interest in these structures is two-fold : they have naturally appeared in Drinfel's approach to a quasi-classical limit of "quantum groups" under the name of Hamilton-Lie groups (~\cite{dri}), later renamed Poisson Lie groups. On the other hand, they appeared in the heart of the Integrable system theory under the name of Sklyanin algebras (~\cite{sky1}, ~\cite{sky2}). It should be noted that almost at the same time, a simpler case of a Sklyanin algebra structure was studied in a purely algebraic context by M. Artin and J. Tate.
In our paper we will concentrate on the Poisson (co)homology of Sklyanin algebras with $4$ generators (which is the original case introduced by Sklyanin ~\cite{sky1}).\\
Generally speaking, Sklyanin algebras belong to a class of "Poisson Structures with regular symplectic leaves". This class was introduced and studied in ~\cite{odru} as the next interesting generalization of "regular" Poisson structures. We mean a class of algebras which are a quasi-classical limit of associative algebras with quadratic relations which are flat deformations of polynomial functions in $\mathbb{C}^n$. It was proven in ~\cite{odru} that these Poisson algebras have polynomial Casimirs (functions which commute with any function on the manifold) and the "dimension" of this algebra is equal to the sum of the degrees of the Casimir generators.\\
Let us describe our basic Poisson algebra in detail.\\
Let
$$q_1=\frac{1}{2}(x_1^2 + x_3^2)+kx_2x_4,$$
$$q_2=\frac{1}{2}(x_1^2 + x_4^2)+kx_1x_3,$$
where $k\in\mathbb C.$\\
We will introduce a Poisson structure $\pi$ on $\mathbb C^4$ or, more generally, on $\mathbb C[x_1, x_2, x_3, x_4]$ (in $\mathbb CP^3$, in the formal series ring $\mathbb C[[x_1, x_2, x_3, x_4]]$,...) by the formula :
$$\{f, g\}_{\pi}:=\frac{df\wedge dg\wedge dq_1\wedge dq_2}{dx_1\wedge dx_2\wedge dx_3\wedge dx_4}.$$
Then the brackets between the coordinate functions are defined by (mod $4$) :
$$\{x_i, x_{i+1}\}=k^2x_ix_{i+1}-x_{i+2}x_{i+3};$$
$$\{x_i, x_{i+2}\}=k(x_{i+3}^2-x_{i+1}^2),$$
$i=1, 2, 3, 4.$\\
This algebra will be denoted by $q_4(\mathcal E)$, where $\mathcal E$ represents an elliptic curve, which parameterizes the algebra (via $k$). We can also think of this curve $\mathcal E$ as a geometric interpretation of the couple $q_1=0$, $q_2=0$ embedded in $\mathcal CP^3$ (as was observed in Sklyanin's initial paper). We will follow Sklyanin's version of this embedding, considering the quadrics $Q_1$, $Q_2$ in $\mathbb C^4$ and the complete intersection of these quadrics :
\begin{equation}\label{e}
Q_1=\frac{1}{2}(x_1^2+x_2^2+x_3^2)
\end{equation}
\begin{equation}\label{f}
Q_2=\frac{1}{2}(x_4^2+J_1x_1^2+J_2x_2^2+J_3x_3^2)
\end{equation}
The algebra $q_4(\mathcal E)$ is given by the formulas :
$$\{x_1, x_i\}=(-1)^idet(\frac{\partial Q_k}{\partial x_l}), \ l\neq1, i; \ k=1, 2; $$
$$\{x_i, x_j\}=(-1)^{i+j}det(\frac{\partial Q_k}{\partial x_l}), \ l\neq i, j; \ k=1, 2. $$
More generally, considering $n-2$ polynomials $Q_i$ in $K^n$ with coordinates $x_i$, $i=1,...,n$, where $K$ is a field of characteristic zero, we can define, for any polynomial $\lambda\in K[x_1,...,x_n]$, a bilinear differential operation :\\
$$\{\cdot ,\cdot\} : K[x_1,...,x_n]\otimes K[x_1,...,x_n]\longrightarrow K[x_1,...,x_n]$$
by the formula
\begin{equation}\label{q}
\{f,g\}=\lambda\frac{df\wedge dg\wedge dQ_1\wedge...\wedge dQ_{n-2}}{dx_1\wedge dx_2\wedge...\wedge dx_n},\ \  \  f,g\in K[x_1,...,x_n]
 \end{equation}
This operation gives a Poisson algebra structure on $K[x_1,...,x_n].$
The polynomials $Q_i, i=1,...,n-2$ are Casimir functions for the brackets (\ref{q}) and any Poisson structure on $K^n$, with $n-2$ generic Casimirs $Q_i$, is written in this form. Every Poisson structure of this form is called a Jacobian Poisson structure (JPS) (~\cite{khi1}, ~\cite{khi2}).\\
The case $n=4$ in $(\ref{q})$ corresponds to the classical generalized Sklyanin quadratic Poisson algebra.\\
In ~\cite{pic}, Anne Pichereau gives the Poisson (co)homology of a Poisson structure given on $K^n$ by the formula (\ref{q}) when $n=3$ and where the Casimir $Q_1$ is a weight homogeneous polynomial  with an isolated singularity in zero.\\
The goal of this paper is to find the (co)homology of a Poisson structure given on $K^n$ by the formula (\ref{q}) when $n=4$ and where the Casimirs $Q_1$ and $Q_2$ are weight homogeneous with an isolated singularity in zero. We use a  method similar to that of Pichereau in ~\cite{pic},  Van den Bergh in ~\cite{van}, or Marconnet in ~\cite{mar}.\\
The results of the paper will be used in our subsequent research of cohomological properties of "quantum" counterparts of the algebra $q_4(\mathcal E)$ known also as Feigin-Odesskii-Sklyanin Elliptic algebras $Q_4(\mathcal E)$. We will apply the Brylinsky spectral sequence
arguments to compute their Hochshild homology (~\cite{pelap}).\\
The paper is organized as follows. We start by introducing some basic notions of Poisson (co)homolo\-gy, the De Rham complex and we introduce some applications and operators that we use to have a simple description of K\"{a}hler differentials and of Poisson homology complexes.
The next part is devoted to homological tools we are going to use to compute the Poisson homology.
 The last part of the paper is devoted to the computation of the Poisson (co)homology of polynomial Poisson structures when the space of Casimir functions is generated by two polynomials which form complete intersection with an isolated singularity. We apply our method to compute the Poisson (co)-homology of the Sklyanin algebra.\\
 {\bf Acknowledgements.} The author is grateful to Jean-Claude Thomas, Yves Felix and Michel Granger for useful comments and discussions. This work is a first part of my thesis prepared at the University of Angers. I would like to take this opportunity to thank my advisors, Vladimir Roubtsov and Bitjong Ndombol, for suggesting to me this interesting problem and for their availability during this projet. This work was partially supported by the Programme SARIMA.
\section{Poisson (co)homology complex}
Consider $(\mathcal A, \pi=\{\, \ \})$ a Poisson algebra. This means an antisymmetric biderivation $\{\cdot, \cdot\} : \mathcal A\times\mathcal A\rightarrow\mathcal A$ such that $(\mathcal A, \{\cdot, \cdot\})$ is a Lie algebra.
\subsection{The de Rham complex and Poisson homology complex}
We recall that the $\mathcal A$-module of K\"ahler differentials of $\mathcal A$ is denoted by $\Omega^1(\mathcal A)$ and the graded $\mathcal A$-module $\Omega^p(\mathcal A) :=\bigwedge^p\Omega^1(\mathcal A)$ is the module of all K\"ahler $p$-differential.
As a vector space, respectively, as an $\mathcal A$-module, $\Omega^p(\mathcal A)$ is generated by elements of the form $FdF_1\wedge...\wedge dF_p$, respectively of the form, $dF_1\wedge...\wedge dF_p$, where $F, F_i\in\mathcal A$, $i=1,...,p$.
We denote by $\Omega^{\bullet}(\mathcal A)=\displaystyle{\oplus_{p\in\mathbb N}}\Omega^p(\mathcal A)$, with the convention that $\Omega^0(\mathcal A)=\mathcal A$, the space of all K\"ahler differentials.\\
The differential $d : \mathcal A\longrightarrow\Omega^1(\mathcal A)$ extends to a graded $K$-linear map $$d :\Omega^\bullet(\mathcal A)\longrightarrow\Omega^{\bullet+1}(\mathcal A)$$
by setting :
$$d(GdF_1\wedge...\wedge dF_p) :=dG\wedge dF_1\wedge...\wedge dF_p$$
for $G,F_1,...,F_p\in\mathcal A$, where $p\in\mathbb{N}$. It is called the de Rham  differential. It is a graded derivation, of degree $1$, of $(\Omega^\bullet(\mathcal A), \wedge)$, such that $d^2=0$. The resulting complex is called the de Rham complex and its cohomology is the de Rham cohomology of $\mathcal A$.\\
The Poisson boundary operator, also called the Brylinsky or Koszul differential and denoted by ${\partial : \Omega^\bullet(\mathcal A)\longrightarrow\Omega^{\bullet-1}(\mathcal A)}$, is given by :
$$\partial_k(F_0dF_1\wedge...\wedge dF_k)=\displaystyle{\sum_{1\leq i\leq k}}(-1)^{i+1}\{F_0, F_i\}dF_1\wedge...\wedge \widehat{dF_i}\wedge ...\wedge dF_k$$
$$+\displaystyle{\sum_{1\leq i<j\leq k}} (-1)^{i+j}F_0d\{F_i,F_j\}\wedge dF_1\wedge...\wedge \widehat{dF_i}\wedge...\widehat{dF_j}\wedge...\wedge dF_k$$\\
where $F_0,...,F_k\in\mathcal A.$\\
One can check, by a direct computation, that $\partial_k$ is well-defined and that it is a boundary operator, $\partial_k\circ\partial_{k+1}=0.$\\
The homology of this complex is called the Poisson homology associated to $(\mathcal A, \pi)$ and is denoted by $H_{\bullet}(\mathcal A, \pi).$
\subsection{The Poisson cohomology complex}
\begin{definition} A skew-symmetric $k$-linear map $P\in\mbox{Hom}_K(\wedge^k\mathcal A,\mathcal A)$ is called a skew-symmetric $k$-derivation of $\mathcal A$ with values in $\mathcal A$ if it is a derivation in each of its arguments.
\end{definition}
The $\mathcal A$-module of skew-symmetric $k$-derivation is denoted by $\mathcal X^k(\mathcal A).$ We define the graded $\mathcal A$-module
$$\mathcal X^\bullet(\mathcal A):=\displaystyle{\bigoplus_{k\in\mathbb N}}\mathcal X^k(\mathcal A)$$
whose elements are called skew-symmetric multi-derivations. By convention, the first term in this sum, $\mathcal X^0$, is $\mathcal A$, and $\mathcal X^k(\mathcal A):=\{0\}$ for $k<0$.\\
The Poisson coboundary operator associated with $(\mathcal A, \pi)$ and denoted by $\delta : \mathcal X^\bullet(\mathcal A)\longrightarrow\mathcal X^{\bullet+1}(\mathcal A)$, is given by :
$$\delta^k(Q)(F_0, F_1,..., F_k)=\displaystyle{\sum_{1\leq i\leq k}}(-1)^{i}\{F_i, Q(F_0, F_1,...,\widehat F_i,..., F_k)\}$$
$$+\displaystyle{\sum_{1\leq i<j\leq k}} (-1)^{i+j}Q(\{F_i, F_j\}, F_0, ...,\widehat F_i,..., \widehat F_j,...,F_k)$$\\
where $F_0,...,F_k\in\mathcal A,$ and $Q\in\mathcal X^k(\mathcal A).$\\
One can check, by a direct computation, that $\delta^k$ is well-defined and that it is a coboundary o\-pe\-ra\-tor, $\delta^{k+1}\circ\delta^{k}=0.$\\
The cohomology of this complex is called the Poisson cohomology associated with $(\mathcal A, \pi)$ and denoted by $H^{\bullet}(\mathcal A, \pi).$
\subsection{Unimodular Poisson structure}
In this part, we consider the affine space of dimension $n$ $K^n$ and its algebra of regular functions $\mathcal A=K[x_1,...,x_n]$.\\
We denote by $S_{p,q}$ the set of all $(p,q)-$shuffles, that is permutations $\sigma$ of the set $\{1,...,p+q\}$ such that $\sigma(1)<...<\sigma(p)$ and $\sigma(p+1)<...<\sigma(p+q), p, q\in\mathbb N $.\\
The family of maps $\star : \mathcal X^k(\mathcal A)\longrightarrow\Omega^{n-k}(\mathcal A)$, defined by :
$$\star Q=\displaystyle{\sum_{\sigma\in S_{k,n-k}}}\epsilon(\sigma)Q(x_{\sigma(1)},...,x_{\sigma(k)})dx_{\sigma(k+1)}\wedge...\wedge dx_{\sigma(n)}$$
are isomorphisms.\\
Assume a Poisson structure on $\mathcal A$ is given. It is natural to ask wether we have the same duality between the Poisson cohomology and the Poisson homology. Generally, the answer to this question is negative. Besides, it is easy to see that the answer depends on the Poisson structure we have. For example such a duality does exist for the wide and important class of unimodular Poisson structure (~\cite{xu1}) which Jacobian Poisson structures belong (~\cite{khi2}).\\
If we denote by $D_\bullet$ the map : $D_\bullet := \star^{-1}\circ d\circ\star : \mathcal X^\bullet(\mathcal A)\longrightarrow\mathcal X^{\bullet -1}(\mathcal A)$, where $d$ is the de Rham differential, we say that a Poisson bracket $\pi$ on $\mathcal A$ is unimodular if $D_2(\pi)=0.$\\
Our purpose is to find the homology and the cohomology of  the Jacobian structures on\\ $K[x_1, x_2, x_3, x_4]$. These structures being unimodular, our objective in the sequel will be to determine the Poisson homology of such Poisson structures.
\subsection{Vector notations}
Now, we are going to present some vector notations which we shall use afterwards. Let us consider the following applications and differential operators : \\ \\

$\begin{array}{cccc}
  \times : & \mathcal A^4\times\mathcal A^4 &\longrightarrow & \mathcal A^6 \\
    &\left(\begin{array}{ccc}
                         \overrightarrow{X}=\left(
                             \begin{array}{c}
                               X_1 \\
                               \cdot\\
                               \cdot \\
                               X_4 \\
                             \end{array}
                           \right)
          & ,&\overrightarrow{Y}=\left(
                             \begin{array}{c}
                               Y_1 \\
                               \cdot\\
                               \cdot \\
                               Y_4 \\
                             \end{array}
                           \right)\end{array}\right)&                    \longmapsto & \overrightarrow{X}\times\overrightarrow{Y}=\left(
                                                                \begin{array}{c}
                                                                  X_1Y_4-X_4Y_1 \\
                                                                X_1Y_2-X_2Y_1 \\
                                                                  X_3Y_2-X_2Y_3 \\
                                                                  X_3Y_4-X_4Y_3\\
                                                                  X_3Y_1-X_1Y_3\\
                                                                  X_2Y_4-X_4Y_2\\
                                                                \end{array}
                                                              \right)
\end{array}$\\ \\ \\
$\begin{array}{cccc}
  \bar{\times} : & \mathcal A^4\times\mathcal A^6 &\longrightarrow & \mathcal A^4 \\
    &\left(\begin{array}{ccc}
                         \overrightarrow{X}=\left(
                             \begin{array}{c}
                               X_1 \\
                               \cdot\\
                               \cdot \\
                               X_4 \\
                             \end{array}
                           \right)
          & ,&\overrightarrow{Y}=\left(
                             \begin{array}{c}
                               Y_1 \\
                               \cdot\\
                               \cdot\\
                               \cdot\\
                               Y_6 \\
                             \end{array}
                           \right)\end{array}\right)&                    \longmapsto & \overrightarrow{X}\bar{\times}\overrightarrow{Y}=\left(
                                                                \begin{array}{c}
                                                                  -X_4Y_3+X_2Y_4-X_3Y_6\\
                                                                  \space X_3Y_1-X_1Y_4+X_4Y_5\\
                                                                  -X_2Y_1+X_4Y_2+X_1Y_6\\
                                                                  -X_3Y_2+X_1Y_3-X_2Y_5\\
                                                                \end{array}
                                                              \right)
\end{array}$\\ \\ \\
We denote by $"\cdot"$ the scalar product in $\mathcal A^4$ or in $\mathcal A^6.$\\

$f : \mathcal A^6\longrightarrow \mathcal A^6$ is  an $\mathcal A$-linear morphism given by the matrix :
$$\left(
  \begin{array}{cccccc}
    0 & 0 & -1 & 0 & 0 & 0 \\
    0& 0 & 0 & 1 & 0 & 0 \\
    -1 & 0 & 0 & 0 & 0 & 0 \\
    0 & 1 & 0 & 0 & 0 & 0 \\
    0 & 0 & 0 & 0 & 0 & 1 \\
    0 & 0 & 0 & 0 & 1 & 0 \\
  \end{array}
\right)$$\\ \\ \\
$\begin{array}{lr}\begin{array}{cccc}
   \overrightarrow{\nabla} :& \mathcal A &\longrightarrow & \mathcal A^4 \\
    &F&                    \longmapsto &\overrightarrow{\nabla}F=\left(
                                \begin{array}{c}
                                  \frac{\partial F}{x_1} \\
                                  \cdot \\
                                  \cdot \\
                                  \frac{\partial F}{x_4} \\
                                \end{array}
                              \right)\
    \end{array}&\begin{array}{cccc}
  \overrightarrow\nabla\times : & \mathcal A^4 &\longrightarrow & \mathcal A^6 \\
    &\overrightarrow{Y}=\left(\begin{array}{c}
                               Y_1 \\
                               \cdot\\
                               \cdot \\
                               Y_4 \\
                             \end{array}
                           \right)
                           &\longmapsto & \overrightarrow\nabla\times\overrightarrow{Y}=\left(
                                                                \begin{array}{c}
                                                                  \frac{\partial Y_4}{\partial x_1}-\frac{\partial Y_1}{\partial x_4} \\
                                                                \frac{\partial Y_2}{\partial x_1}-\frac{\partial Y_1}{\partial x_2} \\
                                                                  \frac{\partial Y_2}{\partial x_3}-\frac{\partial Y_3}{\partial x_2} \\
                                                                  \frac{\partial Y_4}{\partial x_3}-\frac{\partial Y_3}{\partial x_4}\\
                                                                  \frac{\partial Y_1}{\partial x_3}-\frac{\partial Y_3}{\partial x_1}\\
                                                                  \frac{\partial Y_4}{\partial x_2}-\frac{\partial Y_2}{\partial x_4}\\
                                                                \end{array}
                                                              \right)
\end{array}\end{array}$\\ \\ \\
   $\begin{array}{cccc}
   \overrightarrow{\nabla}\bar{\times} :& \mathcal A^6 &\longrightarrow & \mathcal A^4 \\
    &\overrightarrow{G}=\left(
                         \begin{array}{c}
                           G_1\\
                           \cdot \\
                           \cdot\\
                           \cdot \\
                           G_6 \\
                         \end{array}
                       \right)
    &                    \longmapsto &\overrightarrow{\nabla}\bar{\times}\overrightarrow{G}=\left(
                                \begin{array}{c}
                                  -\frac{\partial G_3}{\partial x_4}+\frac{\partial G_4}{\partial x_2}-\frac{\partial G_6}{\partial x_3}\\
                                  \space\frac{\partial G_1}{\partial x_3}-\frac{\partial G_4}{\partial x_1}+\frac{\partial G_5}{\partial x_4} \\
                                  -\frac{\partial G_1}{\partial x_2}+\frac{\partial G_2}{\partial x_4}+\frac{\partial G_6}{\partial x_1}\\
                                  -\frac{\partial G_2}{\partial x_3}+\frac{\partial G_3}{\partial x_1}-\frac{\partial G_5}{\partial x_2}\\
                                \end{array}
                              \right)\
    \end{array}$\\ \\ \\
$\begin{array}{cccc}
   \mbox{Div}(\cdot) :& \mathcal A^4 &\longrightarrow & \mathcal A \\
    &\overrightarrow{K}=\left(
                         \begin{array}{c}
                           K_1\\
                           \cdot \\
                           \cdot\\
                           K_4 \\
                         \end{array}
                       \right)
    &                    \longmapsto &\mbox{Div}(\overrightarrow{K})=\displaystyle{\sum_{i=1}^{4}}\frac{\partial K_i}{\partial x_i}
    \end{array}$\\ \\
    By direct computation, we obtain the following properties :
   \begin{prop}\label{pr1}
The previous operators satisfy the following properties :\\
\begin{enumerate}
  \item[1.]\label{al1} $f^2=Id_{\mathcal A^6};$
  \item[2.]\label{al2} $f(\vec x)\cdot\vec y=\vec x\cdot f(\vec y),$ \ \ $\vec x, \vec y\in\mathcal A^6;$
  \item[3.]\label{al3} $f(\vec x)\cdot f(\vec y)=\vec x\cdot \vec y,$ \ \ $\vec x, \vec y\in\mathcal A^6;$
  \item[4.]\label{al4} $\vec x\cdot(\vec y\times\vec z)=\vec y\cdot(\vec z\bar\times f(\vec x)),$\ \ $\vec x\in\mathcal A^6,$ $\vec y, \vec z\in\mathcal A^4; $
  \item[5.]\label{al5} $\vec x\cdot(\vec y\bar\times\vec z)=-\vec y\cdot(\vec x\bar\times\vec z),$ \ \ $\vec x,\vec y\in\mathcal A^4,$ $\vec z\in\mathcal A^6;$
  \item[6.]\label{al6} $\vec x\cdot(\vec x\bar\times\vec z)=0,$ $\vec x\in\mathcal A^4,$ \ \ $\vec z\in\mathcal A^6;$
  \item[7.]\label{al7} $\vec x\bar\times(\vec x\times\vec y)=0,$ \ \ $\vec x, \vec y\in\mathcal A^4$
  \item[8.]\label{al8} $(\vec x\times\vec z)\cdot f(\vec x\times\vec y)=0,$ \ \ $\vec x, \vec y, \vec z\in\mathcal A^4;$
  \item[9.]\label{al9} $\vec x\bar\times(\vec y\times\vec z)=\vec y\bar\times(\vec z\times\vec x),$ \ \ $\vec x, \vec y, \vec z\in\mathcal A^4;$
  \item[10.]\label{al10} $(\vec x\bar\times f(\vec y\times\vec z))\bar\times(\vec y\times\vec z)=0,$\ \ $\vec x, \vec y, \vec z\in\mathcal A^4;$
  \item[11.]\label{al11} $\vec z\bar\times f(\vec x\times\vec y)=-(\vec z\cdot\vec x)\vec y +(\vec z\cdot\vec y)\vec x,$ \ \ $\vec x, \vec y , \vec z\in\mathcal A^4$;
  \item[12.]\label{al12} $\left(\vec z\bar\times f(\vec x\times\vec y)\right)\times\vec t=-(\vec z\cdot\vec x)\vec y\times\vec t +(\vec z\cdot\vec y)\vec x\times\vec t,$ \ \ $\vec x, \vec y , \vec z, \vec t\in\mathcal A^4$;
  \item[13.]\label{al13} $(\vec x\bar\times f(\vec y\times\vec z))\times\vec z=(\vec x\cdot\vec z)\vec y\times\vec z,$\ \ $\vec x, \vec y, \vec z\in\mathcal A^4;$
      \item[14.]\label{al14} $(\vec x\bar\times\vec z)\bar\times f(\vec x\times\vec y)=-(\vec z\cdot f(\vec x\times\vec y))\vec x,$\ \ $\vec x, \vec y, \vec z\in\mathcal A^4, \vec z\in\mathcal A^6;$
  \item[15.]\label{al15} $\overrightarrow{\nabla}\bar\times(\vec x\times\vec y)=\vec y\bar\times(\overrightarrow{\nabla}\times\vec x)-\vec x\bar\times(\overrightarrow{\nabla}\times\vec y),$\ \
$\vec x, \vec y\in\mathcal A^4;$
  \item[16.]\label{al16} $\overrightarrow{\nabla}\times F\vec x=F\overrightarrow{\nabla}\times\vec x+\overrightarrow{\nabla}F\times\vec x,$\ \ $F\in\mathcal A,$ $\vec x\in\mathcal A^4;$
  \item[17.]\label{al17} $\overrightarrow{\nabla}\bar\times F\vec y=F\overrightarrow{\nabla}\bar\times\vec y+\overrightarrow{\nabla}F\bar\times\vec y,$ \ \ $F\in\mathcal A,$ $\vec y\in\mathcal A^6;$
  \item[18.]\label{al18} $Div(F\vec x)=\overrightarrow{\nabla}F\cdot\vec x+FDiv(\vec x),$\ \ $F\in\mathcal A,$ $\vec x\in\mathcal A^4;$
\item[19.]\label{al19} $Div(\vec x\bar\times\vec y)=\vec y\cdot f(\overrightarrow{\nabla}\times\vec x)-\vec x\cdot(\overrightarrow{\nabla}\bar\times\vec y),$ $ \vec x\in\mathcal A^4,$\ \ $\vec y\in\mathcal A^6.$
\end{enumerate}
\begin{proof}
Let us give a proof for formula $(19)$ for example.\\
Consider $\vec x= (X_1,X_2,X_3,X_4)^t\in \mathcal A^4$ and $\vec y= (Y_1,...,Y_6)^t\in\mathcal A^6.$\\
We have :
$$\begin{array}{rl}
    Div(\vec x\bar\times\vec y)&=\frac{\partial}{\partial x_1}(-X_4Y_3+X_2Y_4-X_3Y_6)+\frac{\partial}{\partial x_2}(X_3Y_1-X_1Y_4+X_4Y_5)\\
    &\ \ + \frac{\partial}{\partial x_3}(-X_2Y_1+X_4Y_2+X_1Y_6)+\frac{\partial}{\partial x_4}(-X_3Y_2+X_1Y_3-X_2Y_5)\\
    & \\
   &=Y_1(\frac{\partial X_3}{\partial x_2}-\frac{\partial X_2}{\partial x_3})+Y_2(\frac{\partial X_4}{\partial x_3}-\frac{\partial X_3}{\partial x_4})+Y_3(\frac{\partial X_1}{\partial x_4}-\frac{\partial X_4}{\partial x_1})\\
   & \ \ +Y_4(\frac{\partial X_2}{\partial x_1}-\frac{\partial X_1}{\partial x_2})+Y_5(\frac{\partial X_4}{\partial x_2}-\frac{\partial X_2}{\partial x_4})+Y_6(\frac{\partial X_1}{\partial x_3}-\frac{\partial X_3}{\partial x_1})\\
   &\ \ -X_1(\frac{\partial Y_4}{\partial x_2}-\frac{\partial Y_6}{\partial x_3}-\frac{\partial Y_3}{\partial x_4})-X_2(-\frac{\partial Y_4}{\partial x_1}+\frac{\partial Y_1}{\partial x_3}+\frac{\partial Y_5}{\partial x_4})\\
   &\ \ -X_3(\frac{\partial Y_6}{\partial x_1}-\frac{\partial Y_1}{\partial x_2}+\frac{\partial Y_2}{\partial x_4})-X_4(\frac{\partial Y_3}{\partial x_1}-\frac{\partial Y_5}{\partial x_2}-\frac{\partial Y_2}{\partial x_3})\\
   & \\
   &=\vec y\cdot f(\overrightarrow{\nabla}\times\vec x)-\vec x\cdot(\overrightarrow{\nabla}\bar\times\vec y)
\end{array}$$

\end{proof}
   \end{prop}
According to the definition of K\"ahler differentials, we have the following isomorphisms of $\mathcal A$-modules \\
$$\begin{array}{ccc}\Omega^1(\mathcal A)&\stackrel{\sim}{\longrightarrow} &\mathcal A^4\\
F_1dx_1+F_2dx_2+F_3dx_3+F_4dx_4&\longmapsto &(F_1,..., F_4)\\
\end{array}$$\\
$$\begin{array}{ccc}\Omega^2(\mathcal A)&\stackrel{\sim}{\longrightarrow} &\mathcal A^6\\
F_1dx_1\wedge dx_4+F_2dx_1\wedge dx_2+F_3dx_3\wedge dx_2&\longmapsto &(F_1, F_2,..., F_6)\\
+F_4dx_3\wedge dx_4+F_5dx_3\wedge dx_1+F_6dx_2\wedge dx_4& &\\
\end{array}$$\\
$$\begin{array}{ccc}\Omega^3(\mathcal A)&\stackrel{\sim}{\longrightarrow} &\mathcal A^4\\
K_1dx_2\wedge dx_3\wedge dx_4+K_2dx_3\wedge dx_1\wedge dx_4&\longmapsto &(K_1,..., K_4)\\
+K_3dx_1\wedge dx_2\wedge dx_4+K_4dx_2\wedge dx_1\wedge dx_3& &\\
\end{array}$$\\
$$\begin{array}{ccc}\Omega^4(\mathcal A)&\stackrel{\sim}{\longrightarrow} &\mathcal A\\
Udx_1\wedge dx_2\wedge dx_3\wedge dx_4&\longmapsto &U\\
\end{array}$$\\
According to the previous isomorphisms, we can write the de Rham complex in terms of elements of $\mathcal A$, $\mathcal A^4$ and $\mathcal A^6$ :\\
\begin{equation}\label{qa}
\begin{array}{ccccccccccccc}
  K&\stackrel{d}{\longrightarrow}&\mathcal A & \stackrel{d}{\longrightarrow} & \mathcal A^4 & \stackrel{d}{\longrightarrow} & \mathcal A^6 & \stackrel{d}{\longrightarrow} & \mathcal A^4 & \stackrel{d}{\longrightarrow} & \mathcal A & \stackrel{d}{\longrightarrow} &0\\
   & &F & \longmapsto & \overrightarrow{\nabla}F &  &  &  &  &  &  & & \\
   & && & \overrightarrow{F} &\longmapsto & \overrightarrow{\nabla}\times\overrightarrow{F} & &  &  & & & \\
  & &   & & & &\overrightarrow{G} &\longmapsto & \overrightarrow{\nabla}\bar\times \overrightarrow{G}&  &  &  &  \\
 & &  &  &  &  &  & &\overrightarrow{K} & \longmapsto & \mbox{Div}(\overrightarrow{K}) &
\end{array}
\end{equation}
\begin{prop}{\textnormal{(Poincaré's lemma)}}
The de Rham complex $(\ref{qa})$ of the polynomial algebra $\mathcal A=K[x_1,..,x_4]$ is an exact one.
\end{prop}

\begin{prop}
According to the previous isomorphisms, Poisson boundary operators associated with the Poisson algebra $(\mathcal A,\pi)$ given by two generic Casimir functions $P_1$ and $P_2$, can be written in a compact form :\\
\begin{equation}\label{hc}
0\longrightarrow\mathcal A\stackrel{\partial_4}{\longrightarrow}\mathcal A^4\stackrel{\partial_3}{\longrightarrow}\mathcal A^6\stackrel{\partial_2}{\longrightarrow}\mathcal A^4\stackrel{\partial_1}{\longrightarrow}\mathcal A
\end{equation}\\
where :\\ \\
$\partial_1(\overrightarrow{H})=(\overrightarrow\nabla\times\overrightarrow H)\cdot f(\overrightarrow\nabla P_1\times\overrightarrow\nabla P_2), \overrightarrow H\in\mathcal A^4$\\ \\
$\partial_2(\overrightarrow{G})= -(\overrightarrow\nabla\bar{\times}\overrightarrow G)\bar{\times} f(\overrightarrow\nabla P_1\times\overrightarrow\nabla P_2)-\overrightarrow\nabla(\overrightarrow G\cdot f(\overrightarrow\nabla P_1\times\overrightarrow\nabla P_2)), \overrightarrow G\in\mathcal A^6$\\ \\
$\partial_3(\overrightarrow K)=\mbox{Div}(\overrightarrow K)\overrightarrow\nabla P_1\times\overrightarrow\nabla P_2+\overrightarrow\nabla\times(\overrightarrow K\bar\times f(\overrightarrow\nabla P_1\times\overrightarrow\nabla P_2))$\\ \\
$\partial_4(U)=-\overrightarrow\nabla U\bar\times(\overrightarrow\nabla P_1\times\overrightarrow\nabla P_2)$
\end{prop}
\begin{proof}
Let us give a proof of the last formula and let $Udx_1\wedge dx_2\wedge dx_3\wedge dx_4$ be an element of $\Omega^4(\mathcal A).$
We have $\partial_4(Udx_1\wedge dx_2\wedge dx_3\wedge dx_4)= \ (I) \ + \ (II), $ where :
$$(I)= \{U, x_1\}dx_2\wedge dx_3\wedge dx_4+\{U, x_2\}dx_3\wedge dx_1\wedge dx_4+\{U, x_3\}dx_1\wedge dx_2\wedge dx_4+\{U, x_4\}dx_2\wedge dx_1\wedge dx_3$$
and
\begin{multline*}
    (II)=-Ud\{x_1, x_2\}\wedge dx_3\wedge dx_4+Ud\{x_1, x_3\}\wedge dx_2\wedge dx_46-d\{x_1, x_4\}\wedge dx_2\wedge dx_3+\\
    -Ud\{x_2, x_3\}\wedge dx_1\wedge dx_4+Ud\{x_2, x_4\}\wedge dx_1\wedge dx_3-Ud\{x_3, x_4\}\wedge dx_1\wedge dx_2.
\end{multline*}
This second term is exactly $(II)=-Ud(dP_1\wedge dP_2)=0$. Then the boundary $\partial_4(Udx_1\wedge dx_2\wedge dx_3\wedge dx_4)$ is equal to the first term $(I)$.\\
But using our identification, we have $\partial_4(U)= (K_1, K_2, K_3, K_4)^t$, where $K_i=\{U, x_i\}.$\\
We can see, by a simply computation, that
$$K_1=\{U, x_1\}:=\frac{dU\wedge dx_1\wedge dP_1\wedge dP_2}{dx_1\wedge dx_2\wedge dx_3\wedge dx_4}$$
is equal to
$$\frac{\partial U}{\partial x_4}(\frac{\partial P_1}{\partial x_3}\frac{\partial P_2}{\partial x_2}-\frac{\partial P_1}{\partial x_2}\frac{\partial P_2}{\partial x_3})-\frac{\partial U}{\partial x_2}(\frac{\partial P_1}{\partial x_3}\frac{\partial P_2}{\partial x_4}-\frac{\partial P_1}{\partial x_4}\frac{\partial P_2}{\partial x_3})+\frac{\partial U}{\partial x_3}(\frac{\partial P_1}{\partial x_2}\frac{\partial P_2}{\partial x_4}-\frac{\partial P_1}{\partial x_4}\frac{\partial P_2}{\partial x_2})$$
which is exactly the first coordinate of $-\overrightarrow\nabla U\bar\times(\overrightarrow\nabla P_1\times\overrightarrow\nabla P_2).$
\end{proof}
The Poisson homology takes the following form :\\ \\
$\begin{array}{lll}
H_0(\mathcal A, \pi)&=&\frac{\mathcal A}{\{(\overrightarrow\nabla\times\overrightarrow H)\cdot f(\overrightarrow\nabla P_1\times\overrightarrow\nabla P_2)| \overrightarrow H\in\mathcal A^4\}}\\&& \\
H_1(\mathcal A, \pi)&=&\frac{\{\overrightarrow H\in\mathcal A^4|(\overrightarrow\nabla\times\overrightarrow H)\cdot f(\overrightarrow\nabla P_1\times\overrightarrow\nabla P_2)=0\}}{\{-(\overrightarrow\nabla\bar\times\overrightarrow G)\bar\times f(\overrightarrow\nabla P_1\times\overrightarrow\nabla P_2)-\overrightarrow\nabla(\overrightarrow G\cdot f(\overrightarrow\nabla P_1\times\overrightarrow\nabla P_2)\}}\\&&\\
H_2(\mathcal A, \pi)&=&\frac{\{\overrightarrow G\in\mathcal A^6|(\overrightarrow\nabla\bar\times\overrightarrow G)\bar\times f(\overrightarrow\nabla P_1\times\overrightarrow\nabla P_2)+\overrightarrow\nabla(\overrightarrow G\cdot f(\overrightarrow\nabla P_1\times\overrightarrow\nabla P_2))=0\}}{\{{Div}(\overrightarrow K) f(\overrightarrow\nabla P_1\times\overrightarrow\nabla P_2)+\overrightarrow\nabla\times[\overrightarrow K\bar\times f(\overrightarrow\nabla P_1\times\overrightarrow\nabla P_2)]\}}\\&&\\
H_3(\mathcal A, \pi)&=&\frac{\{\overrightarrow K\in\mathcal A^4|{Div}(\overrightarrow K) f(\overrightarrow\nabla P_1\times\overrightarrow\nabla P_2)+\overrightarrow\nabla\times[\overrightarrow K\bar\times (\overrightarrow\nabla P_1\times\overrightarrow\nabla P_2)]=0\}}{\{-\overrightarrow\nabla U\bar\times (\overrightarrow\nabla P_1\times\overrightarrow\nabla P_2)\}}\\&&\\
H_4(\mathcal A, \pi)&=&\{U\in\mathcal A|\overrightarrow\nabla U\bar\times (\overrightarrow\nabla P_1\times\overrightarrow\nabla P_2)=0\}\\
\end{array}$
\section{Homological tools}
In this part, we introduce the homological tools that we will need to find the Poisson homology of our Poisson structure. Here, $\mathcal A$ is the polynomial algebra $K[x_1, x_2, x_3, x_4]$.
\subsection{Weight homogeneous skew-symmetric multi-derivations}
The Schouten bracket is a family of maps
$$[\cdot,\cdot]_S : \mathcal X^p(\mathcal A)\times\mathcal X^q(\mathcal A)\longrightarrow\mathcal X^{p+q-1}(\mathcal A),$$
defined by \\
$\begin{array}{c}
[P, Q]_S(F_1,...,F_{p+q-1})=\displaystyle{\sum_{\sigma\in S_{q,p-1}}}\epsilon(\sigma)P(Q(F_{\sigma(1)},...,F_{\sigma(q)}),F_{\sigma(q+1)},...,F_{\sigma(q+p-1)})\\
\end{array}$\\
$\begin{array}{cccccccccccc}
 & & & & & & & & & & &-(-1)^{(p-1)(q-1)}\displaystyle{\sum_{\sigma\in S_{p,q-1}}}\epsilon(\sigma)Q(P(F_{\sigma(1)},...,F_{\sigma(p)}),F_{\sigma(p+1)},...,F_{\sigma(p+q-1)})\\
\end{array}$\\ \\
for $P\in\mathcal X^p(\mathcal A)$, $Q\in\mathcal X^q(\mathcal A)$, and for $F_1,...,F_{p+q-1}\in\mathcal A$
for $p, q\in\mathbb N.$\\
By convention, $S_{p,-1}:=\emptyset$ and $S_{-1,q}:=\emptyset$, for $p, q\in\mathbb N.$
\begin{definition}
Let $\mathcal V\in\mathcal X^1(\mathcal A)$ and $Q\in\mathcal X^q(\mathcal A)$. Then the Lie derivative of $Q$ with respect to $\mathcal V$ is $\mathcal L_{\mathcal V}Q:=[\mathcal V,Q]_S$
\end{definition}
\begin{definition}
A non-zero multi-derivation $P\in\mathcal X^\bullet(\mathcal A)$ is said to be weight homogeneous of degree $r\in\mathbb Z$, if there exists positive integers $\varpi_1, \varpi_2, \varpi_3, \varpi_4\in\mathbb N^{\star}$, the weights of the variables $x_1, x_2, x_3, x_4$, without a common divisor, such that $$\mathcal L_{\vec{e}_\varpi}(P)=rP$$
where $\mathcal L_{\vec{e}_\varpi}$ is a Lie derivative with respect to the Euler derivation
$$\vec e_\varpi:=\varpi_1x_1\frac{\partial}{\partial x_1}+...+\varpi_4x_4\frac{\partial}{\partial x_4}$$
\end{definition}
The degree of a weight homogeneous multi-derivation $P\in\mathcal X^\bullet(\mathcal A)$ is also denoted by $\varpi(P)\in\mathbb Z$
By convention, the zero $k$-derivation is weight homogeneous of degree $-\infty$.\\
The Euler derivation $\vec e_\varpi$ is identified (with the isomorphisms of the first section) to the element $\vec e_\varpi=(\varpi_1x_1,...,\varpi_4x_4)\in\mathcal A^4$. We denote by $|\varpi|$ the sum of the weights $\varpi_1+\varpi_2+\varpi_3+\varpi_4$, so that $|\varpi|=Div(\vec e_\varpi)$.\\
Euler's formula, for a weight homogeneous $F\in\mathcal A$, can be written as $\overrightarrow{\nabla}F\cdot\vec e_\varpi=\varpi(F)F,$ and yields, using (15) of proposition 2.1 : $Div(F\vec e_\varpi)=(\varpi(F)+|\varpi|)F.$

The operator $\star$ allows us to transport the notion of weight homogeneity of skew-symmetric multi-derivations to K\"ahler $p$-differential forms.\\
Fixing weights $\varpi_1, \varpi_2, \varpi_3, \varpi_4\in\mathbb N^\star$, it is clear that $\mathcal A={\bigoplus_{i\in\mathbb N}}\mathcal A_i$, where $\mathcal A_0=K$ and for $i\in\mathbb N^\star$, $\mathcal A_i$ is the $K$-vector space generated by all weight homogeneous polynomials of degree $i$. Denoting by $\Omega^k(\mathcal A)_i$ the $K$-vector space given by $\Omega^k(\mathcal A)_i=\{P\in\Omega^k(\mathcal A) : \varpi(P)=i\}\cup\{0\}$,we have the following isomorphisms : \\
\begin{align}
    \Omega^4(\mathcal A)_i\simeq\mathcal X^0(\mathcal A)_i&\simeq\mathcal A_i\nonumber\\
     \Omega^3(\mathcal A)_i\simeq\mathcal X^1(\mathcal A)_i&\simeq\mathcal A_{i+\varpi_1}\times\mathcal A_{i+\varpi_2}\times\mathcal A_{i+\varpi_3}\times\mathcal A_{i+\varpi_4}\nonumber\\
     \Omega^2(\mathcal A)_i\simeq\mathcal X^2(\mathcal A)_i&\simeq\mathcal A_{i+\varpi_1+\varpi_4}\times\mathcal A_{i+\varpi_1+\varpi_2}\times\mathcal A_{i+\varpi_3+\varpi_2}\times\mathcal A_{i+\varpi_3+\varpi_4}\times\mathcal A_{i+\varpi_3+\varpi_1}\times\mathcal A_{i+\varpi_2+\varpi_4}\nonumber\\
     \Omega^1(\mathcal A)_i\simeq\mathcal X^3(\mathcal A)_i&\simeq\mathcal A_{i+\varpi_2+\varpi_3+\varpi_4}\times\mathcal A_{i+\varpi_3+\varpi_2+\varpi_4}\times\mathcal A_{i+\varpi_1+\varpi_2+\varpi_4}\times\mathcal A_{i+\varpi_2+\varpi_1+\varpi_3}\nonumber\\
 \Omega^0(\mathcal A)_i\simeq\mathcal X^4\mathcal A)_i&\simeq\mathcal A_{i+\varpi_1+\varpi_2+\varpi_3+\varpi_4}\nonumber
\end{align}
\begin{rem}
Each arrow of the complex given by (\ref{qa}) is a weight homogeneous map of degree zero, while each arrow of the complex given by (\ref{hc}) is a weight homogeneous map of degree $\varpi(P_1) +\varpi(P_2),$ if $P_1$ and $P_2$ are weight homogeneous elements of $\mathcal A.$
\end{rem}
\subsection{The Koszul complex-complete intersection with an isolated singularity}
\begin{definition}
A weight homogeneous element $P\in\mathcal A=K[x_1,x_2,x_3,x_4]$ has an isolated singularity if
\begin{equation}\label{singularity}
    \mathcal A_{sing}(P):=K[x_1,x_2,x_3,x_4]/<\frac{\partial P}{\partial x_1},\frac{\partial P}{\partial x_2},\frac{\partial P}{\partial x_3},\frac{\partial P}{\partial x_4}>
\end{equation}
has a finite dimension as a $K$-vector space.
\end{definition}
The dimension of $\mathcal A_{sing}(P)$ is called the Milnor number of the singular point.\\
We shall now give a definition of dimension for rings. For this purpose, note that the length of the chain $P_r\supset P_{r-1}\supset\cdot\cdot\cdot\supset P_0$ involving $r+1$ distinct ideals of a given ring is taken to be $r.$
\begin{definition}
The Krull dimension of a ring $\mathcal R$ is the supremum of the lengths of chains of distinct prime ideals in $\mathcal R.$
\end{definition}
\begin{definition}
Let $\mathcal R$ be an associative and commutative graded $K$-algebra. A system of homogeneous elements $a_1,...,a_d$ in $\mathcal R$, where $d$ is the Krull dimension of $\mathcal R$, is called a homogeneous system of parameters of $\mathcal R$ (h.s.o.p.) if $\mathcal R/<a_1,...,a_d>$ is a finite dimensional $K$-vector space.
\end{definition}
For example, if we consider the $K$-algebra $\mathcal A=K[x_1,...,x_4]$, graded by the weight degree, we have a natural h.s.o.p. given by the system $x_1,x_2,x_3,x_4$.
 \begin{definition}
A sequence $a_1,...,a_n$ in a commutative associative algebra $\mathcal R$ is said to be an $\mathcal R$-regular sequence if $<a_1,...,a_n>\neq\mathcal R$ and $a_i$ is not a zero divisor of $\mathcal R/<a_1,...,a_{i-1}>$ for $i=1,2,...,n$.
\end{definition}
For any regular sequence $a_1,...,a_n$, we can define a Koszul complex which is exact (see Weibel ~\cite{wei}) :\\ \\
$\begin{array}{c}
0\longrightarrow\bigwedge^0(\mathcal R^n)\longrightarrow\cdot\cdot\cdot\longrightarrow\bigwedge^{n-2}(\mathcal R^n)\stackrel{\wedge \omega}{\longrightarrow}\bigwedge^{n-1}(\mathcal R^n)\stackrel{\wedge \omega}{\longrightarrow}\bigwedge^n(\mathcal R^n)\\
\end{array}$\\ \\
where $\omega=\displaystyle{\sum_{i=1}^n}a_ie_i$ and $(e_1, e_2,\cdots, e_n)$ is a basis of an $\mathcal R$-module free $\mathcal R^n.$\\
In our particular case, $\mathcal R=K[x_1,x_2,\cdots,x_n]$, using the identifications $\bigwedge^p(\mathcal R^n)\backsimeq \Omega^p(\mathcal R)$, the Koszul complex associated to the sequence $\frac{\partial P}{\partial x_1},\frac{\partial P}{\partial x_2},\cdots,\frac{\partial P}{\partial x_n}$ ($P\in\mathcal R$) have the following form :

$\begin{array}{c}
0\longrightarrow\mathcal A\stackrel{\wedge dP}{\longrightarrow}\Omega^1(\mathcal A)\longrightarrow\cdot\cdot\cdot\longrightarrow\Omega^{n-2}(\mathcal A)\stackrel{\wedge dP}{\longrightarrow}\Omega^{n-1}(\mathcal A)\stackrel{\wedge dP}{\longrightarrow}\Omega^n(\mathcal A)\\
\end{array}$\\
Using the vector notation for $n=4$, we have the following complex : \\
$$0\longrightarrow\mathcal A\stackrel{\overrightarrow\nabla P}{\longrightarrow}\mathcal A^4\stackrel{\times\overrightarrow\nabla P}{\longrightarrow}\mathcal A^6\stackrel{\overrightarrow\nabla P\bar\times}{\longrightarrow}\mathcal A^4\stackrel{\cdot\overrightarrow\nabla P}{\longrightarrow}\mathcal A$$
\begin{theo}\textnormal{(Cohen-Macaulay)}. Let $\mathcal R$ be a noetherian graded $K$-algebra. If $\mathcal R$ has a h.s.o.p. which is a regular sequence, then any h.s.o.p. in $\mathcal R$ is a regular sequence.
\end{theo}
Thus, for each $P\in\mathcal A=K[x_1,...,x_4]$ which is a weight homogeneous polynomial with an isolated singularity, the sequence $\frac{\partial P}{\partial x_1},\frac{\partial P}{\partial x_2},\frac{\partial P}{\partial x_3},\frac{\partial P}{\partial x_4}$ is regular, the associated Koszul complex is exact.
\begin{definition}
Let $\mathcal R$ be a noetherian commutation ring with unit. The depth, $dpth(I)$, of an ideal $I$ of $\mathcal R$ is the maximal length $q$ of an $\mathcal R$-regular sequence $a_1,\cdots,a_q\in I.$
\end{definition}
Let $M$ be a free $\mathcal R$-module of finite rank $n$, where $\mathcal R$ is a a noetherian commutative ring with unit. We denote by $\bigwedge^p(M)$ the $p$-$th$ exterior product of $M$. By convention $\bigwedge^0(M)=\mathcal R.$\\
Let $\eta_1,\cdots,\eta_k$ be given elements of $M$, and $(e_1,\cdots,e_n)$ be a basis of $M.$
$$\eta_1\wedge\cdots\wedge\eta_k=\displaystyle{\sum_{1\leq i_1<\cdots i_k\leq n}}a_{i_1,\cdots,i_k}e_{i_1}\wedge\cdots\wedge e_{i_k}.$$
We denote by $\mathfrak{A}$ the ideal of $\mathcal R$ generated by the coefficients $a_{i_1,\cdots,i_k}$, $1\leq i_1<\cdots i_k\leq n.$\\
Then we define :
$Z^p:=\left\{\eta\in\bigwedge^p(M) : \eta\wedge\eta_1\wedge\cdots\wedge\eta_k=0\right\}$, $p=0, 1, 2,\cdots$\\
$H^p:=Z^p\diagup\displaystyle{\sum_{i=1}^k}\eta_i\wedge\bigwedge^{p-1}(M)$, $p=0, 1, 2,\cdots$\\
We have the following result from Kyoji Saito :
\begin{theo}\textnormal{(~\cite{sai})}
$H^p=0$ for $0\leq p<dpth(\mathfrak{A}).$
\end{theo}
Let us give an example. Suppose $\mathcal A=K[x_1, x_2, x_3, x_4]$ and consider $P_1, P_2$, two weight homogeneous polynomials in $\mathcal A$. We say that $(P_1, P_2)$ defines a complete intersection if $(P_1, P_2)$ is a regular sequence in $\mathcal A$. And $(P_1, P_2)$ has an isolated singularity if $\mathcal A/\langle P_1, P_2, \frac{\partial P_1}{\partial x_i}\frac{\partial P_2}{\partial x_j}-\frac{\partial P_1}{\partial x_j}\frac{\partial P_2}{\partial x_i}, i<j=1, 2, 3, 4\rangle$ is a finite dimensional $K$-vector space. This dimension is also called the Milnor number of singularity and denoted $\mu$.\\
Let $(P_1, P_2)$ be a complete intersection with an isolated singularity. \\
We denote by $\eta_j=\displaystyle{\sum_{i=1}^4}\frac{\partial P_j}{\partial x_i}e_i$, $j=1, 2$, where $(e_1, e_2, e_3, e_4)$ is a basis of a free  $\mathcal A$-module $\mathcal A^4.$\\
Then $\eta_1\wedge\eta_2=\displaystyle{\sum_{i<j=1}^4}a_{ij}e_i\wedge e_j$, where $a_{ij}=\frac{\partial P_1}{\partial x_i}\frac{\partial P_2}{\partial x_j}-\frac{\partial P_1}{\partial x_j}\frac{\partial P_2}{\partial x_i}.$\\
Let $\mathfrak{A}=(a_{ij}, i<j=1,\cdots, 4).$
From the book of {\footnotesize {E.J.N. LOOIJENGA}} (~\cite{loo}, pages 25 and 49), $dpth(\mathfrak{A})=3$.\\
Then for $\eta\in\bigwedge^p(\mathcal A^4)$, $p=1, 2$, if $\eta\wedge\eta_1\wedge\eta_2=0$, we have $\eta=\alpha_1\wedge\eta_1+\alpha_2\wedge\eta_2$, $\alpha_1, \alpha_2\in\bigwedge^{p-1}(\mathcal A^4)$.
Using the vector notation, we get the following results :
\begin{lem}\label{l1}
For $\overrightarrow{G}\in\mathcal A^6,$ $(\overrightarrow\nabla P_1\times\overrightarrow\nabla P_2)\cdot f(\overrightarrow{G})=0$ if and only if $\overrightarrow{G}=\overrightarrow{H}_1\times\overrightarrow\nabla P_1 + \overrightarrow{H}_2\times\overrightarrow\nabla P_2$, where $\overrightarrow{H}_1, \overrightarrow{H}_2\in\mathcal A^4 .$
\end{lem}
\begin{lem}\label{l2}
For $\overrightarrow{H}\in\mathcal A^4,$ $\overrightarrow{H}\bar\times(\overrightarrow\nabla P_1\times\overrightarrow\nabla P_2)=0$ if and only if $\overrightarrow{H}=U_1\overrightarrow\nabla P_1 + U_2\overrightarrow\nabla P_1$, where $U_1, U_2\in\mathcal A.$
\end{lem}
\section{Poisson homology and complete intersection with an isolated singularity}
Let us consider the polynomial algebra $\mathcal A=K[x_1,...,x_4]$ where $K$ is a field of characteristic $0$ equipped with the Jacobian Poisson structure $\pi$ given by two weight homogeneous polynomials $P_1, P_2$ which define a complete intersection with an isolated singularity.\\
First we compute the kernels $(ker\partial_i)_{i=1, 2, 3, 4}.$
\begin{prop}
$ker\partial_4\simeq K[P_1,P_2]$
\end{prop}
\begin{proof}
Let $U$ be a weight homogeneous element of $\Omega^4(\mathcal A)\simeq\mathcal A$ such that $\partial_4(U)=0$.
Then $\overrightarrow{\nabla}U\bar\times\left(\overrightarrow{\nabla}P_1\times\overrightarrow{\nabla}P_2 \right)=0.$ Using lemma \ref{l2}, $\overrightarrow{\nabla}U=\alpha_1\overrightarrow{\nabla}P_1+\alpha_2\overrightarrow{\nabla}P_2$, $\alpha_1, \alpha_2\in\mathcal A.$
$\overrightarrow{\nabla}\times\overrightarrow{\nabla}U=\overrightarrow{\nabla}\alpha_1\times\overrightarrow{\nabla}P_1+
\overrightarrow{\nabla}\alpha_2\times\overrightarrow{\nabla}P_2.$ Since $\overrightarrow{\nabla}\times\overrightarrow{\nabla}U=0$, we obtain
$\overrightarrow{\nabla}\alpha_1\bar\times(\overrightarrow{\nabla}P_1\times\overrightarrow{\nabla}P_2)=0$ and $\overrightarrow{\nabla}\alpha_2\bar\times(\overrightarrow{\nabla}P_1\times\overrightarrow{\nabla}P_2)=0.$
Thus $\alpha_1, \alpha_2$ are elements of $\ker\partial_4.$\\
We can continue this procedure by defining elements $(U^{(n)}_i)$ of $\ker\partial_4$
in such a way that $\overrightarrow{\nabla}U^{(n)}_i=U^{(n+1)}_{2i-1}\overrightarrow{\nabla}P_1+U^{(n+1)}_{2i}\overrightarrow{\nabla}P_2.$
It is clear that $max(degU^{(n+1)}_{2i-1}, degU^{(n+1)}_{2i})<degU^{(n)}_i.$\\
Therefore there exists $m\in\mathbb N$ such that $\overrightarrow{\nabla}U_i^{(m)}=0$, $\forall i=1,...,2^m.$ Hence $U_i^{(m)}=C^{(m)}_i\in K.$ Since
$\overrightarrow{\nabla}U^{(m-1)}_i=U^{(m)}_{2i-1}\overrightarrow{\nabla}P_1+U^{(m)}_{2i}\overrightarrow{\nabla}P_2=
C^{(m)}_{2i-1}\overrightarrow{\nabla}P_1+C^{(m)}_{2i}\overrightarrow{\nabla}P_2$, by the Poincaré lemma there exists $C_i^{(m-1)}\in K$ such that $U_i^{(m-1)}=C_{2i-1}^{(m)}P_1+C_{2i}^{(m)}P_2+C_i^{(m-1)}.$\\
We have $U_{2i+1}^{(m-1)}=C_{4i+1}^{(m)}P_1+C_{4i+2}^{(m)}P_2+C_{2i+1}^{(m-1)}$; $U_{2i+2}^{(m-1)}=C_{4i+3}^{(m)}P_1+C_{4i+4}^{(m)}P_2+C_{2i+2}^{(m-1)}$ and $\overrightarrow{\nabla}U^{(m-2)}_{i+1}=U^{(m-1)}_{2i+1}\overrightarrow{\nabla}P_1+U^{(m-1)}_{2i+2}\overrightarrow{\nabla}P_2.$ By an easy computation $C_{4i+3}^{(m)}P_1\overrightarrow{\nabla}P_2+C_{4i+2}^{(m)}P_2\overrightarrow{\nabla}P_2=\overrightarrow{\nabla}V$, $V\in\mathcal A.$
Because $\overrightarrow{\nabla}\times[C_{4i+3}^{(m)}P_1\overrightarrow{\nabla}P_2+
C_{4i+2}^{(m)}P_2\overrightarrow{\nabla}P_1]=\overrightarrow{\nabla}\times\overrightarrow{\nabla}V=0$, we have $C_{4i+3}^{(m)}=C_{4i+2}^{(m)},$
and $U_{i+1}^{(m-2)}=\frac{1}{2}C_{4i+1}^{(m)}P_1^2+\frac{1}{2}C_{4i+4}^{(m)}P_2^2+C_{2i+1}^{(m-1)}P_1+C_{2i+2}^{(m-1)}P_2+C_{4i+2}^{(m)}P_1P_2+C_{i+1}^{(m-2)}. $
At the end of this procedure, we obtain the existence of a $g\in K[P_1,P_2]$ such that $U=g(P_1,P_2).$
\end{proof}
\begin{prop}
$ker\partial_1=\{\alpha_1\overrightarrow\nabla P_1+\alpha_2\overrightarrow\nabla P_2+\overrightarrow\nabla\alpha_3| \alpha_1, \alpha_2, \alpha_3\in\mathcal A\}$
\end{prop}
\begin{proof}
Let $\overrightarrow H(H_1,H_2,H_3,H_4)$ be a weight homogeneous element of $\Omega^1(\mathcal A)\simeq\mathcal A^4$ such that $\partial_1(\overrightarrow H)=(\overrightarrow\nabla\times\overrightarrow H)\cdot f(\overrightarrow\nabla P_1\times\overrightarrow\nabla P_2)=0$.
According to lemma \ref{l1}, there exists $\overrightarrow{H}',\overrightarrow{H}''\in\mathcal A^4$ such that ${\overrightarrow{\nabla}\times\overrightarrow{H}=\overrightarrow{H}'\times\overrightarrow{\nabla}P_1+
\overrightarrow{H}''\times\overrightarrow{\nabla}P_2.}$\\
We have $0=\overrightarrow{\nabla}\bar\times\left(\overrightarrow{\nabla}\times\overrightarrow{H}\right)=
\overrightarrow{\nabla}\bar\times\left(\overrightarrow{H}'\times\overrightarrow{\nabla}P_1\right)+
\overrightarrow{\nabla}\bar\times\left(\overrightarrow{H}''\times\overrightarrow{\nabla}P_2\right)$\\
But using property $(15)$ of proposition \ref{pr1}, we obtain $\overrightarrow{\nabla}P_1\bar\times(\overrightarrow{\nabla}\times\overrightarrow{H}')+
\overrightarrow{\nabla}P_2\bar\times(\overrightarrow{\nabla}\times\overrightarrow{H}'')=0.$\\
As a direct consequence, $\overrightarrow{\nabla}P_2\cdot\left(\overrightarrow{\nabla}P_1\bar\times(\overrightarrow{\nabla}\times\overrightarrow{H}')\right)=0=
\overrightarrow{\nabla}P_1\cdot\left(\overrightarrow{\nabla}P_2\bar\times(\overrightarrow{\nabla}\times\overrightarrow{H}'')\right)$. Therefore $\overrightarrow{H}'$ and $\overrightarrow{H}''$ are other elements of $ker\partial_1$.\\
When we apply this procedure to $\overrightarrow{H}'$ and $\overrightarrow{H}''$, we get four other of elements of $ker\partial_1$.\\
Continuing in this way yields the existence of elements $\left(\overrightarrow{H}^{(l)}_r \right)$ of $ker\partial_1$ such that  $\overrightarrow{\nabla}\times\overrightarrow{H}^{(n)}_i=\overrightarrow{H}^{(n+1)}_{2i-1}
\times\overrightarrow{\nabla}P_1+
\overrightarrow{H}^{(n+1)}_{2i}\times\overrightarrow{\nabla}P_2.$ We can notice that $max\left(deg\overrightarrow{H}^{(n+1)}_{2i-1}, deg\overrightarrow{H}^{(n+1)}_{2i} \right)<deg\overrightarrow{H}^{(n)}_{i}.$ Then there exists $m\in\mathbb N$ such that $\overrightarrow{\nabla}\times\overrightarrow{H}^{(m)}_{i}=0$, $\forall i=1,...,2^m.$
So there exists $\varphi_i^{(m)}\in\mathcal A$ such that $\overrightarrow{H}_i^{(m)}=\overrightarrow{\nabla}\varphi_i^{(m)}.$\\
Since $\overrightarrow{\nabla}\times\overrightarrow{H}^{(m-1)}_i=\overrightarrow{\nabla}\times\left[\varphi_{2i-1}^{(m)}\overrightarrow{\nabla}P_1+\varphi_{2i}^{(m)}\overrightarrow{\nabla}P_2 \right]$, by the Poincaré lemma $\overrightarrow{H}_i^{(m-1)}=\varphi_{2i-1}^{(m)}\overrightarrow{\nabla}P_1+\varphi_{2i}^{(m)}\overrightarrow{\nabla}P_2+
\overrightarrow{\nabla}\varphi_{i}^{(m-1)}$, $\varphi_{i}^{(m-1)}\in\mathcal A$.\\
We have :
$\overrightarrow{H}_{2i+1}^{(m-1)}=\varphi_{4i+1}^{(m)}\overrightarrow{\nabla}P_1+\varphi_{4i+2}^{(m)}\overrightarrow{\nabla}P_2+
\overrightarrow{\nabla}\varphi_{2i+1}^{(m-1)};$ $\overrightarrow{H}_{2i+2}^{(m-1)}=\varphi_{4i+3}^{(m)}\overrightarrow{\nabla}P_1+\varphi_{4i+4}^{(m)}\overrightarrow{\nabla}P_2+
\overrightarrow{\nabla}\varphi_{2i+2}^{(m-1)};$ and
$\overrightarrow{\nabla}\times\overrightarrow{H}^{(m-2)}_{i+1}=\left(\varphi_{4i+3}^{(m)}-\varphi_{4i+2}^{(m)} \right)\overrightarrow{\nabla}P_2\times\overrightarrow{\nabla}P_1+\overrightarrow{\nabla}\times\left[ \varphi_{2i+1}^{(m-1)}\overrightarrow{\nabla}P_1+\varphi_{2i+2}^{(m-1)}\overrightarrow{\nabla}P_2 \right].$
It is easy to see that $\overrightarrow{\nabla}\times\overrightarrow{G}=
\alpha\overrightarrow{\nabla}P_1\times\overrightarrow{\nabla}P_2$, where $\overrightarrow{G}=\overrightarrow{H}_{i+1}^{(m-2)}-\varphi_{2i+1}^{(m-1)}\overrightarrow{\nabla}P_1-
\varphi_{2i+2}^{(m-1)}\overrightarrow{\nabla}P_2$ and ${\alpha=\varphi_{4i+3}^{(m)}-\varphi_{4i+2}^{(m)}.}$\\
Since $0=\overrightarrow{\nabla}\bar\times(\overrightarrow{\nabla}\times\overrightarrow{G})=
\overrightarrow{\nabla}\alpha\bar\times(\overrightarrow{\nabla}P_1\times\overrightarrow{\nabla}P_2)$, $\alpha$ is an element of $ker\partial_4.$
Therefore $\overrightarrow{H}_{i+1}^{(m-2)}=(\varphi_{2i+1}^{(m-1)}+\varphi)\overrightarrow{\nabla}P_1+
\varphi_{2i+2}^{(m-1)}\overrightarrow{\nabla}P_2+
\overrightarrow{\nabla}\psi.$\\
At the end, we have : $\overrightarrow{H}=\alpha_1\overrightarrow{\nabla}P_1+\alpha_2\overrightarrow{\nabla}P_2+
\overrightarrow{\nabla}\alpha_3,$ $\alpha_1, \alpha_2, \alpha_3\in\mathcal A.$
\end{proof}
\begin{prop}
$\ker\partial_3=\{\alpha\vec e_\varpi+\overrightarrow\nabla U\bar\times(\overrightarrow\nabla P_1\times\overrightarrow\nabla P_2)| \alpha\in K[P_1,P_2], U\in\mathcal A\}$, if $|\varpi|=2\varpi(P_1)=2\varpi(P_2).$
\end{prop}
\begin{proof}
Let $\overrightarrow{K}\in ker\partial_3,$ be a weight homogeneous element of $\Omega^3(\mathcal A)\simeq\mathcal A^4.$\\
Then $div(\overrightarrow{K})\overrightarrow{\nabla}P_1\times\overrightarrow{\nabla}P_2+\overrightarrow{\nabla}\times\left(\overrightarrow{K}\bar\times f(\overrightarrow{\nabla}P_1\times\overrightarrow{\nabla}P_2)\right)=0.$
But using property $(8)$ of proposition \ref{pr1}, we see that $\left(\overrightarrow{\nabla}\times(\overrightarrow{K}\bar\times f(\overrightarrow{\nabla}P_1\times\overrightarrow{\nabla}P_2))\right)\cdot f(\overrightarrow{\nabla}P_1\times\overrightarrow{\nabla}P_2)=0,$ In other words $\overrightarrow{K}\bar\times f(\overrightarrow{\nabla}P_1\times\overrightarrow{\nabla}P_2)$ is an element of $ker\partial_1.$ Therefore, ${\overrightarrow{K}\bar\times f(\overrightarrow{\nabla}P_1\times\overrightarrow{\nabla}P_2)=\alpha_1\overrightarrow{\nabla}P_1+
\alpha_2\overrightarrow{\nabla}P_2+\overrightarrow{\nabla}\alpha_3},$ where $\alpha_1, \alpha_2, \alpha_3$ are elements of $\mathcal A.$
Therefore $Div(\overrightarrow{K})\overrightarrow{\nabla}P_1\times\overrightarrow{\nabla}P_2+\overrightarrow{\nabla}\alpha_1\times\overrightarrow{\nabla}P_1
+\overrightarrow{\nabla}\alpha_2\times\overrightarrow{\nabla}P_2=0.$ So $\overrightarrow{\nabla}P_i\bar\times\left[Div(\overrightarrow{K})\overrightarrow{\nabla}P_1\times\overrightarrow{\nabla}P_2+\overrightarrow{\nabla}\alpha_1\times\overrightarrow{\nabla}P_1
+\overrightarrow{\nabla}\alpha_2\times\overrightarrow{\nabla}P_2 \right]=0,$ $i=1,2.$\\
According to properties $(9)$ and $(7)$ of proposition \ref{pr1}, we obtain $\overrightarrow{\nabla}\alpha_i\bar\times(\overrightarrow{\nabla}P_1\times\overrightarrow{\nabla}P_2)=0,$ $i=1,2.$
i.e., $\alpha_1,\alpha_2\in K[P_1,P_2].$\\
Thus $Div(\overrightarrow{K})\overrightarrow{\nabla}P_1\times\overrightarrow{\nabla}P_2+
\overrightarrow{\nabla}\alpha_1\times\overrightarrow{\nabla}P_1
+\overrightarrow{\nabla}\alpha_2\times\overrightarrow{\nabla}P_2=0,$ and we obtain $Div(\overrightarrow{K})=\frac{\partial\alpha_1}{\partial P_2}-\frac{\partial\alpha_2}{\partial P_1}.$
On the other hand, using property $(10)$ of proposition \ref{pr1}, we observe that  $\alpha_3$ is an element of $K[P_1,P_2].$\\
We obtain\\
$$\left\{
\begin{array}{lr}
\overrightarrow{K}\bar\times f(\overrightarrow{\nabla}P_1\times\overrightarrow{\nabla}P_2)=\beta_1
\overrightarrow{\nabla}P_1+
\beta_2\overrightarrow{\nabla}P_2&(a)\\
Div(\overrightarrow{K})=\frac{\partial\beta_1}{\partial P_2}-\frac{\partial\beta_2}{\partial P_1}&(b)\\
\end{array}
\right.$$\\
Here $\beta_1=\alpha_1+\frac{\partial\alpha_3}{\partial P_1},$ and $\beta_2=\alpha_2+\frac{\partial\alpha_3}{\partial P_2}.$\\
$(a)\Rightarrow\left(\overrightarrow{K}\bar\times f(\overrightarrow{\nabla}P_1\times\overrightarrow{\nabla}P_2)
\right)\times\overrightarrow{\nabla}P_i=\left(\beta_1\overrightarrow{\nabla}P_1+
\beta_2\overrightarrow{\nabla}P_2\right)\times\overrightarrow{\nabla}P_i,$ $i=1,2.$\\
By property $(12)$ of proposition \ref{pr1}, we have $\beta_j=(-1)^i\overrightarrow{K}\cdot\overrightarrow{\nabla}P_i,$ $i\neq j=1,2.$\\
According to the degree equation $\beta_1\in P_2K[P_1,P_2]$, and $\beta_2\in P_1K[P_1,P_2].$\\
Let us suppose that $\beta_1=cP_1^{r_1}P_2^{r_2+1}.$\\
Then $\vec e_{\varpi}\cdot\overrightarrow{\nabla}P_2=\varpi(P_2)P_2\Rightarrow \overrightarrow{K}=\frac{c}{\varpi(P_2)}P_1^{r_1}P_2^{r_2}\vec e_\varpi+\overrightarrow{G}$, where $\overrightarrow{G}=\overrightarrow{\nabla}P_2\bar\times\overrightarrow{X},$ $\overrightarrow{X}\in\mathcal A^6.$\\
$\overrightarrow{K}\cdot\overrightarrow{\nabla}P_1=-\beta_2\Rightarrow
\beta_2+cP_1^{r_1+1}P_2^{r_2}=-(\overrightarrow{\nabla}P_2
\bar\times\overrightarrow{X})\cdot\overrightarrow{\nabla}P_1$, $\forall r_1,r_2\in\mathbb N.$\\
According to the degree equation, $\beta_2=-cP_1^{r_1+1}P_2^{r_2}.$\\
Then
$$\left\{
\begin{array}{l}
\overrightarrow{K}\cdot\overrightarrow{\nabla}P_2=cP_1^{r_1}P_2^{r_2+1}\\
\overrightarrow{K}\cdot\overrightarrow{\nabla}P_1=cP_1^{r_1+1}P_2^{r_2}\\
Div(\overrightarrow{K})=c(2+r_1+r_2)P_1^{r_1}P_2^{r_2}\\
\end{array}
\right. $$
$$\overrightarrow{K}=\frac{c}{\varpi(P_2)}P_1^{r_1}P_2^{r_2}\vec e_\varpi+\overrightarrow{G}\Rightarrow\left\{
\begin{array}{l}
\overrightarrow{G}\cdot\overrightarrow{\nabla}P_2=0\\
\overrightarrow{G}\cdot\overrightarrow{\nabla}P_1=0\\
Div(\overrightarrow{G})=c\left(2-\frac{|\varpi|}{\varpi(P_1)}\right)P_1^{r_1}P_2^{r_2}=0\\
\end{array}
\right.$$
Then $\overrightarrow{G}\cdot\overrightarrow{\nabla}P_1=0\Rightarrow\overrightarrow{G}=
\overrightarrow{\nabla}
P_1\bar\times\overrightarrow{X}_1,$ $\overrightarrow{X}_1\in\mathcal A^6.$\\
Using lemma \ref{l1}, $\overrightarrow{X}_1
=\overrightarrow{F}_1\times\overrightarrow{\nabla}P_1+\overrightarrow{F}_2
\times\overrightarrow{\nabla}P_2,$ $\overrightarrow{F}_1, \overrightarrow{F}_2\in\mathcal A^4.$
$$\begin{array}{rl}
\overrightarrow{G}&=\overrightarrow{\nabla}P_1\bar\times\left(\overrightarrow{F}_1\times\overrightarrow{\nabla}P_1+\overrightarrow{F}_2
\times\overrightarrow{\nabla}P_2 \right)\\
&=\overrightarrow{\nabla}P_1\bar\times\left(\overrightarrow{F}_2
\times\overrightarrow{\nabla}P_2 \right)\\
&=\overrightarrow{F}\bar\times\left(\overrightarrow{\nabla}P_1
\times\overrightarrow{\nabla}P_2 \right), \mbox{where} \overrightarrow{F}=-\overrightarrow{F}_2\\
\end{array} $$
Since $Div(\overrightarrow{G})=0$, we have $\left(\overrightarrow{\nabla}P_1\times\overrightarrow{\nabla}P_2
\right)\cdot f(\overrightarrow{\nabla}\times\overrightarrow{F})=0.$
Therefore $\overrightarrow{F}=\alpha_1\overrightarrow{\nabla}P_1+\alpha_2\overrightarrow{\nabla}P_2
+\overrightarrow{\nabla}\alpha_3,$ $\alpha_1, \alpha_2, \alpha_3\in\mathcal A,$
and $\overrightarrow{G}=\overrightarrow{\nabla}\alpha_3\bar\times\left(\overrightarrow{\nabla}P_1\times\overrightarrow{\nabla}P_2\right).$
\end{proof}
\begin{prop}
$ker\partial_2=\{U\overrightarrow\nabla P_1\times\overrightarrow\nabla P_2+\overrightarrow\nabla\alpha_1\times\overrightarrow\nabla P_1+\overrightarrow\nabla\alpha_2\times\overrightarrow\nabla P_2| U, \alpha_1, \alpha_2\in\mathcal A\}$
\end{prop}
\begin{proof}
Let $\overrightarrow{G}\in ker\partial_2,$ be a weight homogeneous element of $\Omega^2(\mathcal A)\simeq\mathcal A^6.$\\
We have $\left(\overrightarrow{\nabla}\bar\times\overrightarrow{G}\right)\bar\times f(\overrightarrow{\nabla}P_1\times\overrightarrow{\nabla}P_2)+
\overrightarrow{\nabla}\left(\overrightarrow{G}\cdot f(\overrightarrow{\nabla}P_1\times\overrightarrow{\nabla}P_2)\right)=0.$\\
Then $\overrightarrow{\nabla}\times\left[\left(\overrightarrow{\nabla}\bar\times\overrightarrow{G}\right)\bar\times f(\overrightarrow{\nabla}P_1\times\overrightarrow{\nabla}P_2)\right]+Div(\overrightarrow{\nabla}\bar\times\overrightarrow{G})
\overrightarrow{\nabla}P_1\times\overrightarrow{\nabla}P_2=0.$\\
 Therefore, $\overrightarrow{\nabla}\bar\times\overrightarrow{G}=\alpha\vec e_\varpi+\overrightarrow{\nabla}U\bar\times(\overrightarrow{\nabla}P_1\times\overrightarrow{\nabla}P_2),$ $\alpha\in K[P_1,P_2],$ $U\in\mathcal A.$\\ $0=Div(\overrightarrow{\nabla}\bar\times\overrightarrow{G})=(\varpi(\alpha)+|\varpi|)\alpha\Rightarrow\alpha=0.$\\
Thus $\overrightarrow{\nabla}\bar\times\overrightarrow{G}=
\overrightarrow{\nabla}U\bar\times(\overrightarrow{\nabla}P_1
\times\overrightarrow{\nabla}P_2),$
i.e., $\overrightarrow{G}=U\overrightarrow{\nabla}P_1
\times\overrightarrow{\nabla}P_2+\overrightarrow{\nabla}\times\overrightarrow{X},$ $\overrightarrow{X}\in\mathcal A^4.$\\
But using properties $(10)$ and $(8)$ of proposition \ref{pr1}, we obtain $\overrightarrow{\nabla}\left((\overrightarrow{\nabla}\times\overrightarrow{X})\cdot f(\overrightarrow{\nabla}P_1\times\overrightarrow{\nabla}P_2)\right)=0$, and $(\overrightarrow{\nabla}\times\overrightarrow{X})\cdot f(\overrightarrow{\nabla}P_1\times\overrightarrow{\nabla}P_2)\in K.$\\
For degree reasons, $(\overrightarrow{\nabla}\times\overrightarrow{X})\cdot f(\overrightarrow{\nabla}P_1\times\overrightarrow{\nabla}P_2)=0.$\\
Then from proposition $4.1$, $\overrightarrow{X}=\alpha_1\overrightarrow\nabla P_1+\alpha_2\overrightarrow\nabla P_2+\overrightarrow\nabla\alpha_3,$ $\alpha_1, \alpha_2, \alpha_3\in\mathcal A.$
\end{proof}
\begin{prop}
The following complexes are exact and the arrows are maps of degree zero :
\begin{equation}\label{z}
    0\longrightarrow K[P_1, P_2]\stackrel{\alpha}{\longrightarrow}\mathcal A(-\varpi(P_1))\oplus\mathcal A(-\varpi(P_2))\oplus\mathcal A\stackrel{\beta}{\longrightarrow}ker\partial_1\longrightarrow 0
\end{equation}\\
$\alpha(u(P_1, P_2))=(-\frac{\partial u}{\partial P_1}, -\frac{\partial u}{\partial P_2}, u)$; \\ $\beta(\alpha_1, \alpha_2, \alpha_3)=\alpha_1\overrightarrow\nabla P_1+\alpha_2\overrightarrow\nabla P_2+\overrightarrow\nabla\alpha_3.$ \\
\begin{equation}\label{y}
    0\longrightarrow K[P_1, P_2](-\varpi(P_1))\oplus K[P_1, P_2](-\varpi(P_2))\stackrel{\gamma}{\longrightarrow}
\end{equation}
\begin{equation*}
    \stackrel{\gamma}{\longrightarrow}\mathcal A(-\varpi(P_1)-\varpi(P_2))\oplus\mathcal A(-\varpi(P_1))\oplus\mathcal A(-\varpi(P_2))\stackrel{\epsilon}{\longrightarrow}ker\partial_2\longrightarrow 0
\end{equation*}\\
$\gamma(\alpha_1, \alpha_2)=(\frac{\partial\alpha_1}{\partial P_2}-\frac{\partial\alpha_2}{\partial P_1}, \alpha_1, \alpha_2)$; \\ $\epsilon(U, \alpha_1, \alpha_2)=U\nabla P_1\times\overrightarrow\nabla P_2+\nabla\alpha_1\times\overrightarrow\nabla P_1+\nabla\alpha_2\times\overrightarrow\nabla P_2.$\\
\begin{equation}\label{x}
   0\longrightarrow K[P_1, P_2](-\varpi(P_1)-\varpi(P_2))\stackrel{\theta}{\longrightarrow}K[P_1, P_2](-|\varpi|)\oplus\mathcal A(-\varpi(P_1)-\varpi(P_2))\stackrel{\sigma}{\longrightarrow}ker\partial_3\longrightarrow 0
\end{equation}\\
$\theta(V(P_1, P_2)=(0, V);$\\
$\sigma(U, V)=U\vec e_\varpi+\overrightarrow\nabla V\bar\times(\overrightarrow\nabla P_1\times\overrightarrow\nabla P_2),$
if $|\varpi|=2\varpi(P_1)=2\varpi(P_2),$\\
\end{prop}
\begin{proof}
Let us give the proof for the first sequence :\\
If $\alpha_1\overrightarrow\nabla P_1+\alpha_2\overrightarrow\nabla P_2+\overrightarrow\nabla\alpha_3=0$, then $(\alpha_1\overrightarrow\nabla P_1+\alpha_2\overrightarrow\nabla P_2+\overrightarrow\nabla\alpha_3)\bar\times( \overrightarrow\nabla P_1\times\overrightarrow\nabla P_2)=0.$\\
Therefore, $\overrightarrow\nabla\alpha_3\bar\times( \overrightarrow\nabla P_1\times\overrightarrow\nabla P_2)=0$ and $\alpha_3\in K[P_1, P_2].$\\
Because 0= $(\alpha_1\overrightarrow\nabla P_1+\alpha_2\overrightarrow\nabla P_2+\overrightarrow\nabla\alpha_3)\times\overrightarrow\nabla P_i$, for $i=1, 2,$ we have $\alpha_i=-\frac{\partial\alpha_3}{\partial P_i}$, for $i=1, 2.$\\
We can conclude that the first complex is exact. The proof of the $2nd$ and $3rd$ complexes is similar.
\end{proof}
We also have the following trivial exact sequence of complexes :\\
\begin{equation}\label{w}
    0\longrightarrow ker\partial_{i+1}\longrightarrow\Omega^{i+1}(\mathcal A)\longrightarrow ker\partial_i\longrightarrow H_i(\mathcal A, \pi)\longrightarrow 0
\end{equation}
where $i=0, 1, 2, 3,$ and the arrows are maps of degree zero.\\
\begin{rem}
By using the exactness of the complexes (\ref{z}), (\ref{y}), (\ref{x}), (\ref{w}), we obtain the Poincaré's series of the Poisson homology groups. We can obviously notice that these series do not depend on the polynomials $P_1$ and $P_2$, but on the weights and degrees of $P_1$ and $P_2$.
\end{rem}
We are explicitly going to calculate these Poincaré's series in the quadratic and homogeneous case.
\begin{theo}
If $\varpi_1=...=\varpi_4=1$ and $ \varpi(P_1)=\varpi(P_2)=2$, then as K-vector spaces, $H_i(\Omega^\bullet(\mathcal A), \pi), i=1,2,3,4,$ have the following Poincaré series :\\
$$\begin{array}{lcl}
P(H_0(\mathcal A, \pi),t)&=&\frac{2t^2+4t+1}{(1-t^2)^2} ;\\
&&\\
P(H_1(\mathcal A, \pi),t)&=&\frac{t^4+4t^3+4t^2+4t}{(1-t^2)^2};\\
&&\\
P(H_2(\mathcal A, \pi),t)&=&\frac{2t^4+4t^3}{(1-t^2)^2}; \\
&&\\
P(H_3(\mathcal A, \pi),t)&=&\frac{t^4}{(1-t^2)^2} ;\\
&&\\
P(H_4(\mathcal A, \pi),t)&=&\frac{t^4}{(1-t^2)^2}. \\
\end{array}$$
\end{theo}
\ \newline
Let $\delta=dx_1\wedge dx_2\wedge dx_3\wedge dx_4$, and $\rho=\varpi_1x_1dx_2\wedge dx_3\wedge dx_4+\varpi_2x_2dx_3\wedge dx_1\wedge dx_4+\varpi_3x_3dx_1\wedge dx_2\wedge dx_4+\varpi_4x_4dx_1\wedge dx_2\wedge dx_3.$\\
\begin{theo}
$H_4(\mathcal A, \pi)$ is a rank $1$ free $K[P_1, P_2]$-module generated by $\delta$.
\end{theo}
\begin{theo} If $|\varpi|=2\varpi(P_1)$, then\\
$H_3(\mathcal A, \pi)$ is a rank $1$ free $K[P_1, P_2]$-module generated by $\rho$.
\end{theo}
\begin{proof}
We have proved that, if $|\varpi|=2\varpi(P_1)$, then $ker\partial_3=Im\partial_4 + K[P_1, P_2]\vec e_\varpi.$\\
Now, let $\alpha\in K[P_1, P_2]$ and $U\in\mathcal A$ such that $\alpha\vec e_\varpi=\overrightarrow\nabla U\bar\times(\overrightarrow\nabla P_1\times\overrightarrow\nabla P_2)$. We can suppose that $\alpha$ is a weight homogeneous polynomial.\\
We have $0=Div(\overrightarrow\nabla U\bar\times(\overrightarrow\nabla P_1\times\overrightarrow\nabla P_2))=Div(\alpha\vec e_\varpi)=(\varpi(\alpha)+|\varpi|)\alpha$.\\
Therefore $\alpha=0.$
\end{proof}
Now we suppose that $P_1= Q_1+Q_2$, $P_2=Q_2$ where :
$Q_1$ and $Q_2$ are given by (\ref{e}) and (\ref{f}). $\varpi_1=...=\varpi_4=1;$ $J_i\neq J_j$ if $i\neq j$; for all $i=1, 2, 3, \ J_i\neq -1, 0, 1.$
\begin{theo}~\cite{boo}
 The homological group $H_0(\mathcal A, \pi)$ is a rank 7 free $K[P_1, P_2]$-module generated by $(\mu_i)_{0\leq i\leq 6} = (1, x_1, x_2, x_3, x_4, x_1^2,x_3^2).$
\end{theo}
\begin{rem}
We can remark that $H_0(\mathcal A, \pi)\simeq K[P_1, P_2]\otimes\mathcal A_{sing}(P_1, P_2).$\\
Here $\mathcal A_{sing}(P_1, P_2)=\frac{\mathcal A}{J}$, where $J$ is the ideal generated by $P_1, P_2$ and the $2\times 2$ minors of the Jacobian matrix of $(P_1, P_2).$
\end{rem}
\begin{theo}
The homological group $H_1(\mathcal A, \pi)$ is a free $K[P_1, P_2]$ module given by :
$$H_1(\mathcal A, \pi) \simeq (\displaystyle{\bigoplus_{k=1}^6} K[P_1, P_2]\overrightarrow{\nabla}\mu_k)\oplus(\displaystyle{\bigoplus_{k=1}^5} K[P_1, P_2]\mu_k\overrightarrow{\nabla}P_1)\oplus K[P_1, P_2]\overrightarrow{\nabla}P_1\oplus K[P_1, P_2]\overrightarrow{\nabla}P_2$$
where, $(\mu_i)_{0\leq i\leq 6} = (1, x_1, x_2, x_3, x_4, x_1^2,x_3^2).$
\end{theo}
\begin{proof}
Let $\overrightarrow{H}\in ker\partial_1.$\\
Then $\overrightarrow{H}=\alpha_1\overrightarrow\nabla P_1+\alpha_2\overrightarrow\nabla P_2+\overrightarrow\nabla\alpha_3,$ $\alpha_i\mathcal A.$\\
Using the previous result :\\
$$\alpha_l= f(\overrightarrow\nabla\times\overrightarrow{H}_l)\cdot(\overrightarrow\nabla P_1\times\overrightarrow\nabla P_2) + \displaystyle{\sum_{i=0}^{N_l}}\displaystyle{\sum_{j=0}^{M_l}}\displaystyle{\sum_{k=0}^{6}}\lambda^l_{i,j,k}P_1^iP_2^j\mu_k$$
For $l=1,2,$ we have :\\
$\begin{array}{ll}
\alpha_l\overrightarrow\nabla P_l&= f(\overrightarrow\nabla\times\overrightarrow{H}_l)\cdot(\overrightarrow\nabla P_1\times\overrightarrow\nabla P_2)\overrightarrow\nabla P_l + \displaystyle{\sum_{i=0}^{N_l}}\displaystyle{\sum_{j=0}^{M_l}}\displaystyle{\sum_{k=0}^{6}}\lambda^l_{i,j,k}P_1^iP_2^j\mu_k\overrightarrow\nabla P_l\\
&=\partial_2\left( (-1)^{l+1}\overrightarrow{H}_l\times\overrightarrow\nabla P_l\right)+\displaystyle{\sum_{i=0}^{N_l}}\displaystyle{\sum_{j=0}^{M_l}}\displaystyle{\sum_{k=0}^{6}}\lambda^l_{i,j,k}P_1^iP_2^j\mu_k\overrightarrow\nabla P_l
\end{array}$\\
$\begin{array}{ll}
\overrightarrow\nabla\alpha_3&= \overrightarrow\nabla\left(f(\overrightarrow\nabla\times\overrightarrow{H}_l)\cdot(\overrightarrow\nabla P_1\times\overrightarrow\nabla P_2)\right) + \displaystyle{\sum_{k=1}^{6}}\left(\displaystyle{\sum_{i=0}^{N_3}}\displaystyle{\sum_{j=0}^{M_3}}\lambda^l_{i,j,k}P_1^iP_2^j\mu_k\right)\overrightarrow\nabla\mu_k\\
& \ \ \ \ \ +\displaystyle{\sum_{i=0}^{N_3}}\displaystyle{\sum_{j=1}^{M_3}}\displaystyle{\sum_{k=0}^{6}}j\lambda^3_{i,j,k}P_1^iP_2^{j-1}\mu_k\overrightarrow\nabla P_2\\
& \ \ \ \ \ \ \ +\displaystyle{\sum_{i=1}^{N_3}}\displaystyle{\sum_{j=0}^{M_3}}\displaystyle{\sum_{k=0}^{6}}i\lambda^3_{i,j,k}P_1^{i-1}P_2^{j}\mu_k\overrightarrow\nabla P_1\\
&=\partial_2(-\overrightarrow\nabla\times\overrightarrow{H}_3)+\displaystyle{\sum_{k=1}^{6}}\left(\displaystyle{\sum_{i=0}^{N_3}}\displaystyle{\sum_{j=0}^{M_3}}\lambda^l_{i,j,k}P_1^iP_2^j\mu_k\right)\overrightarrow\nabla\mu_k\\
&\ \ +\displaystyle{\sum_{k=0}^{6}}\left(\displaystyle{\sum_{i=1}^{N_3}}\displaystyle{\sum_{j=0}^{M_3}}i\lambda^3_{i,j,k}P_1^{i-1}P_2^{j}\mu_k\overrightarrow\nabla P_1 +\displaystyle{\sum_{i=0}^{N_3}}\displaystyle{\sum_{j=1}^{M_3}}j\lambda^3_{i,j,k}P_1^iP_2^{j-1}\mu_k\overrightarrow\nabla P_2\right)
\end{array}$\\
Therefore : $ker\partial_1=Im\partial_2+I_1$, where :\\
$$I_1=\left[ (\displaystyle{\bigoplus_{k=1}^6} K[P_1, P_2]\overrightarrow{\nabla}\mu_k)+(\displaystyle{\bigoplus_{k=0}^6} K[P_1, P_2]\mu_k\overrightarrow{\nabla}P_1)+(\displaystyle{\bigoplus_{k=0}^6} K[P_1, P_2]\mu_k\overrightarrow{\nabla}P_2)\right]$$
Then $H_1(\mathcal A, \pi)= \frac{I_1}{I_1\cap Im\partial_2}.$\\
Let $\overrightarrow{G}_i=f(\overrightarrow{\nabla}x_i\times\vec e_{\varpi})$, $\overrightarrow{G}'_i=f(\overrightarrow{\nabla}x_i^2\times\vec e_{\varpi}),$ for $i=1, 2, 3, 4.$\\
We have $\partial_2(\overrightarrow{G}_i)=-J_ix_i\overrightarrow{\nabla}P_1+(J_i+1)x_i\overrightarrow{\nabla}P_2+ 2\left(J_iP_1-(J_i+1)P_2\right)\overrightarrow{\nabla}x_i,$ for $i=1, 2, 3$ and
$\partial_2(\overrightarrow{G}_4)=-x_4\overrightarrow{\nabla}P_1+x_4\overrightarrow{\nabla}P_2+ 2\left(P_1-P_2\right)\overrightarrow{\nabla}x_4.$\\
Therefore $x_i\overrightarrow{\nabla}P_2$ can be written as a $K[P_1, P_2]$-linear sum of $x_i\overrightarrow{\nabla}P_1$ and $\overrightarrow{\nabla}x_i.$\\
On the other hand,\\
$\partial_2(\overrightarrow{G}_i')=3J_ix_i^2\overrightarrow{\nabla}P_1-(J_i+1)x_i^2\overrightarrow{\nabla}P_2
+4(P_2\overrightarrow{\nabla}P_1-P_1\overrightarrow{\nabla}P_2)+ 4\left(J_iP_1-(J_i+1)P_2\right)\overrightarrow{\nabla}x_i^2,$ for $i=1, 2, 3$ and
$\partial_2(\overrightarrow{G}_4')=3x_4^2\overrightarrow{\nabla}P_1-x_4^2\overrightarrow{\nabla}P_2
+4(P_2\overrightarrow{\nabla}P_1-P_1\overrightarrow{\nabla}P_2)+ 4\left(P_1-P_2\right)\overrightarrow{\nabla}x_4^2.$\\
But using the Casimir relations, $x_4^2= -2J_2P_1 + (J_2+1)P_2 + J_{21}x_1^2 + J_{23}x_4^2$, where $J{ij}=J_i-J_j,$
we obtain :
\begin{multline*}
    \partial_2(\overrightarrow{G}'_4)-\partial_2(\frac{J_{21}}{J_1+1}\overrightarrow{G}'_1)-\partial_2(\frac{J_{23}}{J_3+1}\overrightarrow{G}'_3)=
    \frac{3J_{21}}{J_1+1}x_1^2\overrightarrow{\nabla}P_1+\frac{3J_{23}}{J_3+1}x_3^2\overrightarrow{\nabla}P_1-6J_2P_1\overrightarrow{\nabla}P_1
    +\\+6(J_2+1)P_2\overrightarrow{\nabla}P_1+6J_2P_1\overrightarrow{\nabla}P_2
    -6(J_2+1)P_2\overrightarrow{\nabla}P_2-\frac{4J_{21}}{J_1+1}(P_2\overrightarrow{\nabla}P_1-P_1\overrightarrow{\nabla}P_2)+\\
    -\frac{4J_{21}}{J_1+1}(J_1-(J_1+1)P_2)\overrightarrow{\nabla}x_3^2-
    \frac{4J_{23}}{J_3+1}(P_2\overrightarrow{\nabla}P_1-P_1\overrightarrow{\nabla}P_2)-\frac{4J_{23}}{J_3+1}(J_3-(J_3+1)P_2)\overrightarrow{\nabla}x_3^2
    +\\+4(P_2\overrightarrow{\nabla}P_1-P_1\overrightarrow{\nabla}P_2)-8J_2\left(P_1-P_2\right)\overrightarrow{\nabla}P_1
    +8(J_2+1)\left(P_1-P_2\right)\overrightarrow{\nabla}P_2+\\+4J_{21}\left(P_1-P_2\right)\overrightarrow{\nabla}x_1^2
    +4J_{23}\left(P_1-P_2\right)\overrightarrow{\nabla}x_3^2.
\end{multline*}
Therefore $ker\partial_1=Im\partial_2+I_1'$, where
$$I_1'=(\displaystyle{\sum_{k=1}^6} K[P_1, P_2]\overrightarrow{\nabla}\mu_k)+(\displaystyle{\sum_{k=1}^5} K[P_1, P_2]\mu_k\overrightarrow{\nabla}P_1)+ K[P_1, P_2]\overrightarrow{\nabla}P_1+ K[P_1, P_2]\overrightarrow{\nabla}P_2.$$
Then $H_1(\mathcal A, \pi)=\frac{I_1'}{I_1'\cap Im\partial_2}.$\\
But $H_1(\mathcal A, \pi)$ and $I_1'$ have the same Poincaré series. Thus $P(I_1'\cap Im\partial_2, t)=0$ and $I_2'\cap Im\partial_2=0.$
\end{proof}

\begin{theo} The homological group $H_2(\mathcal A, \pi)$ is a free $K[P_1, P_2]$ module given by :
$$H_2(\mathcal A, \pi) \simeq \left(\displaystyle{\bigoplus_{k=1}^5} K[P_1, P_2](\overrightarrow{\nabla}\mu_k\times\overrightarrow{\nabla}P_1)\right)\oplus K[P_1, P_2](\overrightarrow{\nabla}P_1\times\overrightarrow{\nabla}P_1).$$
where, $(\mu_i)_{1\leq i\leq 6} = (x_1, x_2, x_3, x_4, x_1^2,x_3^2).$
\end{theo}
\begin{proof}
Let $\overrightarrow{G}\in ker\partial_2$.\\
$\overrightarrow{G}=\alpha_0\overrightarrow\nabla P_1\times\overrightarrow\nabla P_2+\overrightarrow\nabla \alpha_1\times\overrightarrow\nabla P_1+\overrightarrow\nabla \alpha_1\times\overrightarrow\nabla P_2,$ $\alpha_0, \alpha_1, \alpha_2\in\mathcal A.$\\
We have
$$\alpha_l= f(\overrightarrow\nabla\times\overrightarrow{H}_l)\cdot(\overrightarrow\nabla P_1\times\overrightarrow\nabla P_2) + \displaystyle{\sum_{i=0}^{N_l}}\displaystyle{\sum_{j=0}^{M_l}}\displaystyle{\sum_{k=0}^{6}}\lambda^l_{i,j,k}P_1^iP_2^j\mu_k$$
But :\\
$\left(f(\overrightarrow\nabla\times\overrightarrow{H}_0)\cdot(\overrightarrow\nabla P_1\times\overrightarrow\nabla P_2)\right)\overrightarrow\nabla P_1\times\overrightarrow\nabla P_2=\partial_3(\overrightarrow{H}_0\bar\times(\overrightarrow\nabla P_1\times\overrightarrow\nabla P_2));$\\
$\partial_3(\overrightarrow\nabla P_1\bar\times(\overrightarrow\nabla\times\overrightarrow{H}_l))= -
\overrightarrow\nabla\left(f(\overrightarrow\nabla\times\overrightarrow{H}_l)\cdot(\overrightarrow\nabla P_1\times\overrightarrow\nabla P_2)\right)\times\overrightarrow\nabla P_l,$ for $l=1, 2.$\\
On the other hand
\begin{multline*}
 \partial_3(\lambda^l_{i,j,k}P_1^iP_2^j\mu_k\vec e_\varpi)=\varpi(\mu_k)\lambda^l_{i,j,k}P_1^iP_2^j\mu_k\overrightarrow\nabla P_1\times\overrightarrow\nabla P_2+\varpi(P_2)\lambda^l_{i,j,k}P_1^iP_2^{j+1}\mu_k\overrightarrow\nabla\mu_k\times\overrightarrow\nabla P_2+\\
 -\varpi(P_1)\lambda^l_{i,j,k}P_1^{i+1}P_2^j\mu_k\overrightarrow\nabla\mu_k\times\overrightarrow\nabla P_2.
\end{multline*}\\
Therefore $ker\partial_2=Im\partial_3 + I_2,$ where
$$I_2= \displaystyle{\sum_{k=1}^6}K[P_1, P_2]\overrightarrow\nabla\mu_k\times\overrightarrow\nabla P_1 + \displaystyle{\sum_{k=1}^6}K[P_1, P_2]\overrightarrow\nabla\mu_k\times\overrightarrow\nabla P_2 + K[P_1, P_2]\overrightarrow\nabla P_1\times\overrightarrow\nabla P_2.$$
Then $H_2(\mathcal A, \pi)=\frac{I_2}{I_2\cap Im\partial_3}.$\\
Let $\overrightarrow K_i=(K_{i1}, K_{i2}, K_{i3}, K_{i4})^t$, where $K_{ij}=\delta_{ij}x_i$, and $\overrightarrow K'_i=(K'_{i1}, K'_{i2}, K'_{i3}, K'_{i4})^t$, where $K'_{ij}=\delta_{ij}.$\\
We have $\partial_3(\overrightarrow K_i')=J_i\overrightarrow{\nabla}x_i\times\overrightarrow{\nabla}P_1 - (J_i+1)\overrightarrow{\nabla}x_i\times\overrightarrow{\nabla}P_2.$\\
Then for $i=1, 2, 3$, $\overrightarrow{\nabla}x_i\times\overrightarrow{\nabla}P_2$ is equal to $\frac{J_i}{J_i+1}\overrightarrow{\nabla}x_i\times\overrightarrow{\nabla}P_1$ modulo $Im\partial_3.$\\
We also have  $\overrightarrow{\nabla}x_4\times\overrightarrow{\nabla}P_2=\partial(-\overrightarrow{K}_4')+\overrightarrow{\nabla}x_4\times\overrightarrow{\nabla}P_1.$\\
On the other hand, $\partial_3(\overrightarrow K_i)=\overrightarrow{\nabla}P_1\times\overrightarrow{\nabla}P_2+J_i\overrightarrow{\nabla}x_i^2\times\overrightarrow{\nabla}P_1 - (J_i+1)\overrightarrow{\nabla}x_i^2\times\overrightarrow{\nabla}P_2$, for $i=1, 2, 3,$ and $\partial_3(\overrightarrow K_4)=\overrightarrow{\nabla}P_1\times\overrightarrow{\nabla}P_2+\overrightarrow{\nabla}x_4^2\times\overrightarrow{\nabla}P_1 - \overrightarrow{\nabla}x_4^2\times\overrightarrow{\nabla}P_2.$\\
But using the Casimir relations, $x_4^2= -2J_2P_1 + (J_2+1)P_2 + J_{21}x_1^2 + J_{23}x_4^2$, where $J{ij}=J_i-J_j,$
we obtain :
\begin{multline*}
    \partial_3(\overrightarrow K_4)=\partial_3(\frac{J_{21}}{J_1+1}\overrightarrow K_1+\frac{J_{23}}{J_3+1}\overrightarrow K_3)-(1+\frac{J_{21}}{J_1+1}+\frac{J_{23}}{J_3+1})\overrightarrow{\nabla}P_1\times\overrightarrow{\nabla}P_2+\\
    +\frac{J_{21}}{J_1+1}\overrightarrow{\nabla}x_1^2\times\overrightarrow{\nabla}P_1+
    \frac{J_{23}}{J_3+1}\overrightarrow{\nabla}x_3^2\times\overrightarrow{\nabla}P_1.
\end{multline*}
Then we can explain $\overrightarrow{\nabla}x_3^2\times\overrightarrow{\nabla}P_1$ as a $K$-linear sum of $\overrightarrow{\nabla}x_1^2\times\overrightarrow{\nabla}P_1$ and $\overrightarrow{\nabla}P_1\times\overrightarrow{\nabla}P_2$ modulo $Im\partial_3.$\\
Therefore $ker\partial_2=Im\partial_3 + I_2',$ where
$$I_2'= \displaystyle{\sum_{k=1}^6}K[P_1, P_2]\overrightarrow\nabla\mu_k\times\overrightarrow\nabla P_1.$$
Then $H_2(\mathcal A, \pi)=\frac{I_2'}{I_2'\cap Im\partial_3}.$\\
But $H_2(\mathcal A, \pi)$ and $I_2'$ have the same Poincaré series. Thus $P(I_2'\cap Im\partial_3, t)=0$ and $I_2'\cap Im\partial_3=0.$
\end{proof}

\

\noindent
Laboratoire Angevin de Recherche en Math\'ematiques Universit\'e D'Angers D\'epartement de Math\'ematiques\\
\it{E-mail address} : \verb"pelap@math.univ-angers.fr "

\end{document}